\documentclass[a4paper]{article}
\usepackage[english]{babel}
\usepackage[utf8]{inputenc}
\usepackage{amsthm}
\usepackage{amsmath} 
\usepackage{amssymb}
\usepackage{mathabx}
\usepackage{mathtools}
\usepackage{caption}
\usepackage{adjustbox}
\usepackage{graphicx}
\usepackage[percent]{overpic}
\usepackage{fancyhdr}
\usepackage{tikz}
\usetikzlibrary{shapes}
\usetikzlibrary{plotmarks}
\usepackage{float}
\usepackage{subcaption}
\usepackage{setspace} 
\usepackage{titlesec}
\usepackage{booktabs}
\usepackage{tabularx}
\usepackage{lastpage}
\usepackage[ruled,linesnumbered]{algorithm2e}
\usepackage{todonotes}
\usepackage{enumitem}
\usepackage{rotating}

\let\chapter\undefined

\titleformat*{\subsection}{\normalsize}
\titleformat*{\subsubsection}{\normalsize}

\newcommand{\mytitle}{Tighter McCormick Relaxations through Subgradient Propagation}
\newcommand{\myshorttitle}{Tighter McCormick Relaxations through Subgradient Propagation}
\newcommand{\myauthor}{Jaromi{\l} Najman $\cdot$ Alexander Mitsos$^{*}$} 
\newcommand{\myauthorshort}{J. Najman, A. Mitsos}
\author{\myauthor}
\usepackage{hyperref} 
\hypersetup{
    bookmarks=true,         
    unicode=false,          
    pdftoolbar=true,        
    pdfmenubar=true,        
    pdffitwindow=false,     
    pdftitle={Tighter McCormick Relaxations through Subgradient Propagation},    
    pdfauthor={Jaromil Najman and Alexander Mitsos},     
    pdfsubject={},   
    pdfcreator={Jaromil Najman},   
    pdfproducer={},  
    pdfkeywords={}, 
    pdfnewwindow=true,      
    colorlinks=false,       
    linkcolor=red,          
    citecolor=green,        
    filecolor=magenta,      
    urlcolor=cyan           
}

\definecolor{rwth}{rgb}{0,0.32,0.62}
\definecolor{rwth-75}{rgb}{0.25,0.49,0.71}
\definecolor{rwth-50}{rgb}{0.55,0.73,0.89}
\definecolor{grun}{rgb}{0.34,0.67,0.15}
\definecolor{rot}{rgb}{0.8,0.02,0.11}
\definecolor{magenta}{RGB}{227,0,102}
\definecolor{petrol}{RGB}{0,97,101}
\definecolor{violett}{RGB}{97,33,88}
\definecolor{maigrun}{RGB}{189 205   0}

{
	\fancyhead[L]{\footnotesize{\myshorttitle}}
	\fancyhead[R]{\footnotesize{\textit{\the\day.\the\month.\the\year}}}
	\fancyfoot[L]{\copyright \footnotesize{\textit{\myauthorshort}}}
	\fancyfoot[R]{\footnotesize{Page \thepage\ of \pageref*{LastPage}}}
	\cfoot{}
}

\fancypagestyle{firststyle}
{
	 \fancyhead[L]{\copyright \footnotesize{\textit{\myauthorshort}}}
	 \fancyhead[C]{}
   \fancyhead[R]{\footnotesize{Page \thepage\ of \pageref*{LastPage}}}
	 \fancyfoot[L]{\footnotesize \rule{0.25\textwidth}{0.4pt} \newline Corresponding author: $^*$A. Mitsos E-mail: amitsos@alum.mit.edu}
   \fancyfoot[C]{}
	 \fancyfoot[R]{}
}

\addto\captionsenglish{%
}

\renewcommand{\vec}[1]{\mathbf{#1}}

\newtheoremstyle{myExamplesRemarks}
  {0pt}   
  {0pt}   
  {\rm}      
  {0pt}       
  {\itshape} 
  {}         
  { }  
{\thmname{#1} \thmnumber{#2} }

\theoremstyle{myExamplesRemarks}
\newtheorem{example}{Example}

\makeatletter
\let\@addpunct\@gobble
\g@addto@macro{\thm@space@setup}{\thm@headpunct{}} 
\makeatother

\begin{document}
\thispagestyle{firststyle}

\begin{flushleft}\begin{large}\textbf{\mytitle}\end{large} \end{flushleft}
\begin{small}\myauthor\end{small}

\begin{flushleft}\begin{small}
RWTH Aachen University\\
AVT - Aachener Verfahrenstechnik\\
Process Systems Engineering\\
Forckenbeckstra{\ss}e 51 \\
D - 52074 Aachen
\end{small}
\end{flushleft}
\vspace{1cm}
\textbf{Abstract} Tight convex and concave relaxations are of high importance in the field of deterministic global optimization. We present a heuristic to tighten relaxations obtained by the McCormick technique. We use the McCormick subgradient propagation (Mitsos et al., SIAM J. Optim., 2009) to construct simple affine under- and overestimators of each factor of the original factorable function. Then, we minimize and maximize these affine relaxations in order to obtain possibly improved range bounds for every factor resulting in possibly tighter final McCormick relaxations. We discuss the heuristic and its limitations, in particular the lack of guarantee for improvement. Subsequently, we provide numerical results for benchmark cases found in the COCONUT library and case studies presented in previous works and discuss computational efficiency. We see that the presented heuristic provides a significant improvement in tightness and decrease in computational time in many cases.  

\section{Introduction}

Methods based on branch-and-bound (B\&B) \cite{Morrison2016} are the state-of-the art algorithms in the field of deterministic global optimization. In general, B\&B methods rely on favorable convergence order of the underlying convex and concave relaxations of all functions involved in a given global optimization problem in order to avoid the so-called cluster effect \cite{du1994cluster,Kannan2017,Kannan2017bab,Kearfott1993,wechsung2014}. Domain and range reduction techniques are employed within B\&B algorithms in order to further increase the quality of the underlying relaxations. These techniques are not necessary to guarantee convergence of the B\&B algorithm, however, they are able to drastically speed up convergence. 

Relaxation techniques based on interval arithmetics \cite{Moore:1979:MAI:1098639,nla.cat-vn1029192} describe a general way to obtain valid over- and underestimations of a multivariate function $f: Z\to \mathbb{R}$. It is a well-known fact that the simplest interval-based method, called natural interval extensions, often provides very loose estimators. Thus, efforts have been made to improve tightness by affine reformulations \cite{comba1993ne,de2004affine,ninin2010reliable} and improved Taylor models \cite{convergenceTaylor2013,sahlodin2011Taylor}.

The method of constructing valid convex and concave relaxations of a continuous factorable function, given by a finite recursion of addition, multiplication and composition via propagation of valid factors, e.g., $F_1\circ f_1 + F_2 \circ f_2 \cdot F_3 \circ f_3$, was presented by McCormick \cite{mccormick1976,mccormick1983nonlinear} and extended to multivariate outer functions $F_i$ in \cite{tsoukalas2014multivariate}. McCormick's idea was used in the development of the well-known auxiliary variable method (AVM) \cite{smith1997global,tawarmalani2002convexification,tawarmalani2005} used in state-of-the-art global optimization solvers such as BARON \cite{tawarmalani2005} and ANTIGONE \cite{misener2014}. 

In order to further improve the tightness of the relaxations constructed with the AVM, many bound-tightening procedures are used such as Optimality-Based Bound Tightening \cite{locatelli2013global}, where additional optimization problems are solved in order to tighten variable bounds; bound propagation techniques \cite{Brearley1975}, where information on a constraint is used to possibly tighten the bounds of a different constraint and finally the variables involved in both constraints; probing \cite{tawarmalani2002convexification}, where valid constraints are derived from non-active constraints and more. The recent article by Puranik and Sahinidis \cite{puranik2017domain} provides a thorough overview of the field of tightening techniques for AVM. Most of the techniques applicable to the AVM are (at least theoretically) directly applicable to the relaxations obtained via the McCormick technique. Still, there are almost no algorithms developed directly for the improvement of relaxations obtained by the McCormick method. Recently, Wechsung et al. \cite{Wechsung2015} present an algorithm for constraint propagation using McCormick relaxations resulting in a reduced variable domain and tighter final McCormick relaxations. They reverse the operations starting with pre-computed McCormick relaxations of a given factorable function $g$ and traverse the factors of $g$ backwards in order to tighten the set of feasible points. Herein, we present a different idea with the same goal of improving the final resulting McCormick relaxations. The presented algorithm uses possibly improved bounds for each factor of $g$ obtained through exploitation of subgradient propagation described by Mitsos et al. in \cite{mitsos2009mccormick} providing a tighter final relaxation. 

The remainder of the manuscript is structured as follows. In Section \ref{sect:def}, we provide basic definitions and notation used throughout the article. We present the algorithm in Section \ref{sect:heuristic} supported by an example and discuss its limitations. Subsequently, we present numerical results in Section \ref{sect:numeric} and examine different adjustments of the presented method. Section \ref{sect:conclusion} concludes the work.

\section{Basic definitions}\label{sect:def}

In the following, if not stated otherwise, we consider a continuous function $f:Z\to \mathbb{R}$ with $Z \in \mathbb{I}\mathbb{R}^n$, where $\mathbb{I}\mathbb{R}$ denotes the set of closed bounded intervals of $\mathbb{R}$. $Z\in \mathbb{I}\mathbb{R}^n$, also called \textit{box}, is defined as $Z\equiv[\vec{z}^L,\vec{z}^U]=[z^L_1,z^U_1]\times \dots \times [z^L_n,z^U_n]$ with $\vec{z}^L,\vec{z}^U \in \mathbb{R}^n$ where the superscripts $L$ and $U$ always denote a lower and upper bound, respectively. We denote the image of $f$ over $Z$ by $f(Z)\in\mathbb{I}\mathbb{R}$. We denote the estimation of the range bounds of $f$ on $Z$ with the use of natural interval extensions by $I_{f,nat}\supset f(Z)$ and the exact bounds by $I_{f,e}=f(Z)$. 

We call a convex function $f^{cv}:Z\to\mathbb{R}$ a convex relaxation (or convex underestimator) of $f$ on $Z$ if $f^{cv}(\vec{z})\leq f(\vec{z})$ for every $\vec{z}\in Z$. Similarly, we call a concave function $f^{cc}:Z\to\mathbb{R}$ a concave relaxation (or concave overestimator) of $f$ on $Z$ if $f^{cc}(\vec{z})\geq f(\vec{z})$ for every $\vec{z}\in Z$. We call the tightest convex and concave relaxations of $f$ the convex and concave envelopes $f^{cv}_e, f^{cc}_e$ of $f$ on $Z$, respectively, i.e., it holds $f^{cv}(\vec{z})\leq f^{cv}_e(\vec{z})\leq f(\vec{z})$ and $f(\vec{z})\leq f^{cc}_e(\vec{z})\leq f^{cc}(\vec{z})$ for all $\vec{z}\in Z$ and all convex relaxations $f^{cv}$ and concave relaxations $f^{cc}$ of $f$ on $Z$, respectively.

For a convex and concave function $f^{cv},f^{cc}:Z\to\mathbb{R}$, we call $\vec{s}^{cv}(\bar{\vec{z}}),\vec{s}^{cc}(\bar{\vec{z}})\in \mathbb{R}^n$ a convex and a concave subgradient of $f^{cv},f^{cc}$ at $\bar{\vec{z}}\in Z$, respectively, if
\begin{alignat}{3}
f^{cv}(\vec{z}) &\geq f^{cv}(\bar{\vec{z}}) + (\vec{s}^{cv}(\bar{\vec{z}}))^T(\vec{z}-\bar{\vec{z}}), &~~\forall \vec{z}\in Z, \tag{A1}\label{aff1} \\
f^{cc}(\vec{z}) &\leq f^{cc}(\bar{\vec{z}}) + (\vec{s}^{cc}(\bar{\vec{z}}))^T(\vec{z}-\bar{\vec{z}}), &~~\forall \vec{z}\in Z, \tag{A2}\label{aff2} 
\end{alignat} 
respectively. We denote the affine functions on the right-hand side of inequalities \eqref{aff1}, \eqref{aff2} constructed with the convex and concave subgradient $\vec{s}^{cv}(\bar{\vec{z}}),\vec{s}^{cc}(\bar{\vec{z}})$ by $f^{cv,sub}(\bar{\vec{z}},\vec{z})$ and $f^{cc,sub}(\bar{\vec{z}},\vec{z})$, respectively. Note that $f^{cv,sub}$ and $f^{cc,sub}$ are valid under- and overestimators of $f$ on $Z$, respectively \cite{mitsos2009mccormick}.

\subsection{McCormick relaxations and subgradient propagation}

We will make use of McCormick propagation rules originally developed by McCormick \cite{mccormick1976} and extended to multivariate compositions of functions by Tsoukalas and Mitsos \cite{tsoukalas2014multivariate}. Rules for the propagation of subgradients for the unviariate composition McCormick rule are given in \cite{mitsos2009mccormick} and of multivariate composition rule in \cite{tsoukalas2014multivariate}. The rules for the construction of McCormick relaxations of binary sum, binary product and univariate composition can be found in, e.g., Propositions 2, 3 and Theorem 5 in \cite{bompadre2012convergence}, Propositions 2.6, 2.7 and Theorem 2.8 in \cite{mitsos2009mccormick} or Section 3 in \cite{mccormick1976}. The corresponding subgradient propagation rules can be found in Proposition 2.9 and Theorems 3.2 and 3.3 in \cite{mitsos2009mccormick}. The rule for the construction of McCormick relaxations with multivariate outer functions can be found in Theorem 2 in \cite{tsoukalas2014multivariate} and the corresponding subgradient propagation rule in Theorem 4 of \cite{tsoukalas2014multivariate}.

\subsection{Computational Graph}\label{ssec:DAG}
\begin{figure}[t]
\centering
\begin{tikzpicture}
\node[shape = rectangle,draw=black, minimum size = 0.75cm,line width=1pt,label={[xshift = 1.3cm,yshift=-1.85cm]$\begin{matrix*}[l] f_1(z)=z \\ \textcolor{rot}{\begin{matrix}f_1^L,f^{cv}_1,s_1^{cv} \\ f_1^U,f^{cc}_1,s_1^{cc}\end{matrix}} \end{matrix*}$}] (z) at (0,0) {$z$};

\node[shape = circle,draw = black, minimum size = 1cm,line width=1pt,label={[xshift = 1.4cm,yshift=-0.85cm]$\begin{matrix*}[l] \textcolor{rot}{\begin{matrix}f_2^L,f^{cv}_2,s_2^{cv} \\ f_2^U,f^{cc}_2,s_2^{cc}\end{matrix}} \\f_2(z)=z^2 \end{matrix*}$}] (sqr) at (-3.5,2) {\^{}$2$};

\node[shape = circle,draw = black, minimum size = 1cm,line width=1pt,label={[xshift = 1.4cm,yshift=-1.7cm]$\begin{matrix*}[l] f_3(z)=z^3 \\ \textcolor{rot}{\begin{matrix}f_3^L,f^{cv}_3,s_3^{cv} \\ f_3^U,f^{cc}_3,s_3^{cc}\end{matrix}}\end{matrix*}$}] (pow3) at (0,2) {\^{}$3$};

\node[shape = circle,draw = black, minimum size = 1cm,line width=1pt,label={[xshift = 1.7cm,yshift=-1.7cm]$\begin{matrix*}[l] f_4(z)=\exp(z) \\ \textcolor{rot}{\begin{matrix}f_4^L,f^{cv}_4,s_4^{cv} \\ f_4^U,f^{cc}_4,s_4^{cc}\end{matrix}} \end{matrix*}$}] (exp) at (3.5,2) {$\exp$};

\node[shape = circle,draw = black, minimum size = 1cm,line width=1pt,label={[xshift = 0.7cm,yshift=-0.1cm]$\begin{matrix*}[l]f_5(z)=z-z^2 \\ \textcolor{rot}{\begin{matrix}f_5^L,f^{cv}_5,s_5^{cv} \\ f_5^U,f^{cc}_5,s_5^{cc}\end{matrix}}\end{matrix*}$}] (minus1) at (-4.5,4) {$-$};

\node[shape = circle,draw = black, minimum size = 1cm,line width=1pt,label={[xshift = 1.5cm,yshift=-1cm]$\begin{matrix}f_6(z)=z^3-\exp(z)\\ \textcolor{rot}{\begin{matrix}f_6^L,f^{cv}_6,s_6^{cv} \\ f_6^U,f^{cc}_6,s_6^{cc}\end{matrix}}\end{matrix}$}] (minus2) at (3.5,4) {$-$};

\node[shape = circle,draw = black, minimum size = 1cm,line width=1pt,label={[xshift = 3cm,yshift=-1cm]$\begin{matrix*}[l]f_7(z)=(z-z^2)\cdot (z^3-\exp(z)) \\ \textcolor{rot}{\begin{matrix}f_7^L,f^{cv}_7,s_7^{cv} \\ f_7^U,f^{cc}_7,s_7^{cc}\end{matrix}}\end{matrix*}$}] (times) at (0,6) {$\times$};

\node[shape = rectangle,draw = black, minimum size = 1cm,line width=1pt] (g) at (0,8) {$g(z)$};

\draw[->,line width=1pt] (z)--(sqr);
\draw[->,line width=1pt] (z)--(pow3);
\draw[->,line width=1pt] (z)-|(exp) ;
\draw[->,line width=1pt] (z)-|(minus1) ;
\draw[->,line width=1pt] (sqr)--(minus1) ;
\draw[->,line width=1pt] (pow3)--(minus2) ;
\draw[->,line width=1pt] (exp)--(minus2);
\draw[->,line width=1pt] (minus1)--(times);
\draw[->,line width=1pt] (minus2)--(times) ;
\draw[->,line width=1pt] (times)--(g) ;
\end{tikzpicture}
\caption{Computational graph for $g(z)=(z-z^2)\cdot(z^3-\exp(z))$ on $Z=[-0.5,1]$.}\label{Fig:DAG}
\end{figure}
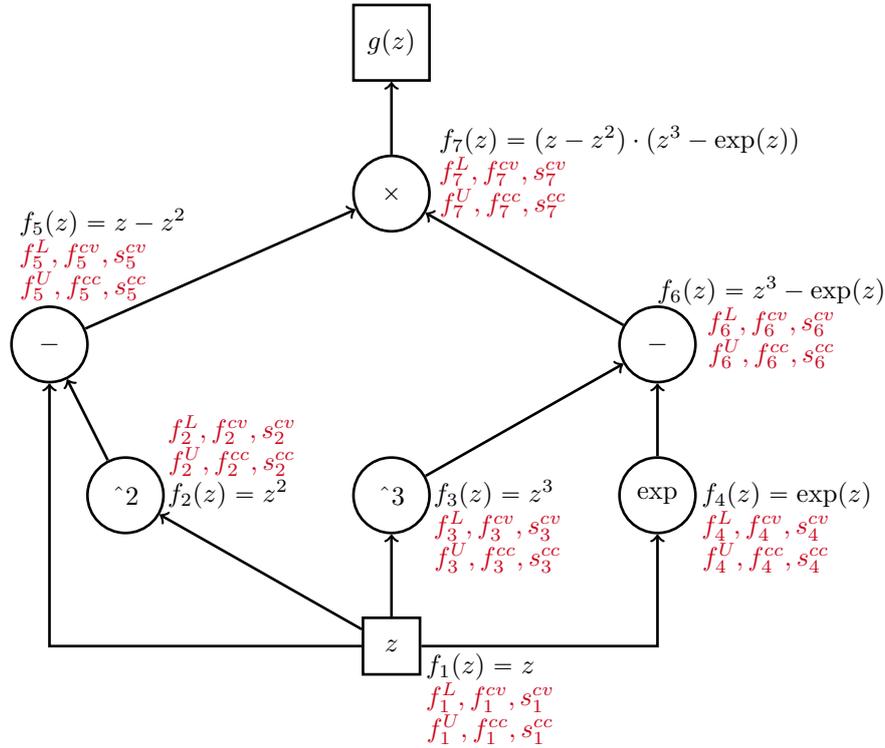

We assume that a directed acyclic graph (DAG) representation $G=(\mathbb{F},E)$, described in, e.g., Sections 2 and 3 in \cite{schichl2005interval}, of a (multivariate) factorable function $g:Z\to \mathbb{R}$ with $Z\in\mathbb{I}\mathbb{R}^n$ is given. $\mathbb{F}$ is the set of vertices, which we call factors herein, consisting of operations and independent variables, and $E$ the set of edges connecting the factors. We assume that for each factor $f_j\in\mathbb{F},j\in\{1,\dots,|\mathbb{F}|\}$ convex $f_j^{cv}(\bar{\vec{z}})$ and concave $f_j^{cc}(\bar{\vec{z}})$ relaxations, the corresponding convex $\vec{s}_j^{cv}(\bar{\vec{z}})$ and concave $\vec{s}_j^{cc}(\bar{\vec{z}})$ subgradients at a point $\bar{\vec{z}} \in Z$ and valid upper $f_j^U$ and lower $f_j^L$ bounds on the range of $f$ over $Z$ are calculated, see Example 4.2 and Fig.~4.2 in \cite{mitsos2009mccormick}. The relaxations and subgradients are calculated by the McCormick rules. The upper $f_j^U$ and lower $f_j^L$ bounds on the range of $f$ over $Z$ are obtained via natural interval extension (\cite{Moore:1979:MAI:1098639,nla.cat-vn1029192}) throughout this article, i.e., $I_{j,nat}=[f_j^L,f_j^U]$. In order to evaluate $g$, its relaxations and its subgradients at a point $\bar{\vec{z}} \in Z$ through $G$, we assume that the corresponding DAG is traversed in a reversed-level-order, i.e., starting at the independent variables $z_i,~i\in\{1,\dots,n \}$ and working through all factors up to the root given as $g(\bar{\vec{z}})$.

\begin{example}
Consider the function $g(z)=(z-z^2)(z^3-\exp(z))$ on $Z=[-0.5,1]$. It consists of 7 factors, namely the independent variable $z$ and the $6$ operations \^{}$2,$\^{}$3,\exp,-,-$ and $\times$. The corresponding computational graph is shown in Fig.~\ref{Fig:DAG}. 
\end{example}

\section{Heuristic for tighter McCormick relaxations}\label{sect:heuristic}

\subsection{Basic idea}

For a nonlinear factorable function given by a finite recursion of addition, multiplication and composition, $g=f_1 \circ f_2+f_3 \circ f_4 \cdot f_5 \circ f_6$, there are several bound and domain tightening techniques and ideas, found in, e.g., \cite{brearley1975analysis,cornelius1984computing,hamed1993calculation,hansen1991analytical,puranik2017bounds} just to name a few and more can be found in the recent article by Puranik and Sahinidis \cite{puranik2017domain}. Many tightening methods use information on constraints within a given problem in order to tighten variable bounds, e.g., \cite{brearley1975analysis,hamed1993calculation,Shectman1998}, while other methods use optimality conditions, reduced costs of variables, and dual multipliers of constraints of the given problem to obtain a tighter relaxation, e.g., \cite{puranik2017bounds,ryoo1996branch,tawarmalani2002convexification}. We present an algorithm which uses information on McCormick relaxations and subgradients of each factor of a particular function $g$ within a given optimization problem to possibly improve the resulting final McCormick relaxations of $g$. In Section 2.3 of \cite{Najman2017anchoring}, we have presented that using tighter range bounds for each factor of a McCormick relaxation results in tighter relaxations. Herein, we present an idea for obtaining tighter McCormick relaxations with the use of subgradients for McCormick relaxations, \cite{mitsos2009mccormick} (implemented within MC\texttt{++}(v2.0)\cite{Chachuatmc++}). The presented algorithm is not guaranteed to improve the final McCormick relaxations making it a heuristic. We first give the basic idea followed by an example and then formalize the algorithm.

Typically, when calculating McCormick relaxations, natural interval extensions are used for the computation of valid interval bounds for the range of each factor of $g$. This is also the case within MC\texttt{++}(v2.0)\cite{Chachuatmc++}. It is possible to use different interval arithmetics, e.g., the standard centered form or Taylor forms (Sections 2.2 and 3.7 in \cite{nla.cat-vn1029192}). We show that the presented algorithm can still improve the final McCormick relaxations  for more sophisticated interval arithmetics and briefly discuss results in Section \ref{ssec:otherIntervals}. By solving $\min\limits_{\vec{z} \in Z} f_j^{cv}(\vec{z})$ and $\max\limits_{\vec{z} \in Z} f_j^{cc}(\vec{z})$ for each factor $f_j, j\in \{1,\dots, |\mathbb{F}|\}$ of a factorable function $g$, where $f_j^{cv},f_j^{cc}$ are McCormick relaxations of $f_j$, we can obtain valid and possibly tighter range bounds for each factor $f_j$. We can achieve even tighter bounds by solving $\min\limits_{\vec{z} \in Z} f_j(\vec{z})$ and $\max\limits_{\vec{z} \in Z} f_j(\vec{z})$ resulting in a possible improvement of McCormick relaxations of $g$ (see Example \ref{ex:heuristic} and Example 5 in \cite{Najman2017anchoring}). However, the number of factors in a factorable function can be very large leading to a high computational time. Thus, we want to find a good trade-off between tightness of range bounds of each factor and computational time needed. We could approximately solve $\min\limits_{\vec{z} \in Z} f_j(\vec{z})$ and $\max\limits_{\vec{z} \in Z} f_j(\vec{z})$ using linear or higher order approximations of $f_j$ in order to simplify the optimization problem but this does not guarantee valid bounds. We could as well approximately solve $\min\limits_{\vec{z} \in Z} f_j^{cv}(\vec{z})$ and $\max\limits_{\vec{z} \in Z} f_j^{cc}(\vec{z})$ by the use of a solution method for convex (nonsmooth) problems, e.g., bundle-methods (\cite{bertsekas2015convex,bertsekas2003convex}), and allow only a small number of iteration steps providing valid but possibly extremely loose bounds. When computing McCormick relaxations of a factorable function $g$, each factor has information about its range bounds, convex and concave McCormick relaxations and its convex and concave subgradients. In this article, we use this information in order to solve $\min\limits_{\vec{z} \in Z} f_j^{cv}(\vec{z})$ and $\max\limits_{\vec{z} \in Z} f_j^{cc}(\vec{z})$ approximately by solving the simple linear box-constrained problems $\min\limits_{\vec{z} \in Z} f_j^{cv,sub}(\bar{\vec{z}},\vec{z})$ and $\max\limits_{\vec{z} \in Z} f_j^{cc,sub}(\bar{\vec{z}},\vec{z})$ for every factor $f_j,j\in \{1,\dots, |\mathbb{F}|\}$ resulting in possibly improved range bounds $f_j^L,f_j^U$ and finally tighter convex and concave McCormick relaxations of $g$. 

In Example 4.4 of \cite{mitsos2009mccormick}, simple affine relaxations are constructed by using the propagated subgradient at a point $\bar{\vec{z}}\in Z$ to construct the affine relaxations $g^{cv,sub}(\bar{\vec{z}},\vec{z}),g^{cc,sub}(\bar{\vec{z}},\vec{z})$ of the original function $g$. The lower bound obtained by evaluating the affine functions at their minimum and maximum, respectively, can result in tighter bounds than the underlying (natural) interval extensions, see Fig.~4.4 in \cite{mitsos2009mccormick}. We can exploit this property of the affine estimators when constructing McCormick relaxations by computing the subgradients at a point in the domain, e.g., the middle-point, in each factor $f_j$ of the factorable function $g$ and checking if we can improve the corresponding range bounds for the current factor, i.e., we approximately solve $\min\limits_{\vec{z} \in Z} f_j^{cv}(\vec{z})$ and $\max\limits_{\vec{z} \in Z} f_j^{cc}(\vec{z})$ by solving $\min\limits_{\vec{z}\in Z} f_j^{cv,sub}(\bar{\vec{z}},\vec{z})$ and $\max\limits_{\vec{z}\in Z} f_j^{cc,sub}(\bar{\vec{z}},\vec{z})$ for every $j\in \{1,\dots, |\mathbb{F}|\}$. We obtain the corresponding values $f_{j,alg}^L,f_{j,alg}^U$ and check if $f_{i,new}^L>f_j^L,f_{j,alg}^U< f_j^U$, i.e., if we can improve the bounds on the range of $f_j$ for each $j\in \{1,\dots, |\mathbb{F}|\}$. Note that it is not guaranteed that the bounds $f^L_{j,alg},f_{j,alg}^U$ are better than the bounds $f_j^L, f_j^U$ obtained through natural interval extensions. Subsequently, we can compute a next linearization point and repeat the computations if desired. Note that since we cannot provide any guarantee on the new bounds $f^L_{j,alg},f_{j,alg}^U$, it can take many re-computations in order to achieve an improvement for the range bounds $f_j^L, f_j^U$ of a factor $f_j$ and thus, a maximal number of iterations should be predefined.

The idea can be intuitively described as a specific application of subgradient bundle methods (\cite{bertsekas2015convex,bertsekas2003convex}). We use subgradients to approximate the possibly non-smooth convex/concave relaxations of a given function $g$ representing our \emph{bundle}. We then construct a linear approximation of $g$ and check if the linearization provides a better range bound than the natural interval extension. If a maximum number of iterations is not reached, we use the subgradient information in order to determine a new point. It is also be possible to describe this idea as a modified Sandwich algorithm (4.2 in \cite{tawarmalani2002convexification}) where polyhedral approximations are constructed for convex functions. However, the Sandwich algorithm does not work with propagated subgradients but rather with the differentials of particular convex functions. Moreover, we try to avoid computing subgradients at all corners, in contrast to what is done in Figures 4.1-4.8 in \cite{tawarmalani2002convexification}, since the dimension of $\vec{z}$ can be too large to make this algorithm applicable in each factor of $g$. The idea of the heuristic is presented in the next example.


\subsection{Illustrative examples}

\subsubsection{Natural interval extensions}

\begin{example}\label{ex:heuristic}
Consider again the function $g(z)=(z-z^2)(z^3-\exp(-z))$ on $Z=[-0.5,1]$, Fig.~\ref{Fig:alg:interval2}, and consider the particular three factors $f_5(z)=z-z^2,~f_6(z)=z^3-\exp(z)$ and $f_7(z)=f_5(z)\cdot f_6(z)$. For factors $f_1,f_2,f_3,f_4$, envelopes are known and natural interval extension provide exact range bounds. The convex and concave McCormick relaxations provide envelopes for $f_5$ on $Z$ given as
\begin{align*}
f_5^{cv}(z)=0.5z-0.5 \text{ and } f_5^{cc}(z)= z- z^2.
\end{align*}
The natural interval extensions are not exact for the range of $f_5$ providing $I_{5,nat}=[-1.5,1]$ while the exact range is given as $I_{5,e}=[-0.75,0.25]$. The convex and concave subgradient of $f_5^{cv}$ and $f_5^{cc}$, respectively, at the middle point $0.25$ of $Z$ are $s^{cv}_5(0.25)=0.5$ and $s^{cc}_5(0.25)=0.5$, respectively. We construct the corresponding affine functions 
\begin{align*}
f_5^{cv,sub}(0.25,z)&=f_5^{cv}(0.25)+s^{cv}_5(0.25)(z-0.25) \text{ and }\\
 f_5^{cc,sub}(0.25,z)&=f_5^{cc}(0.25)+s^{cc}_5(0.25)(z-0.25).
\end{align*}
Next, we evaluate the affine functions at their respective minimum and maximum in order to obtain 
\begin{align}
\begin{split}
\min_{z\in [-0.5,1]} f_5^{cv,sub}(0.25,z) &= f_5^{cv,sub}(s^{cv}_5(0.25),-0.5) = -0.75 \text{ and } \\
\max_{z\in [-0.5,1]} f_5^{cc,sub}(0.25,z) &= f_5^{cc,sub}(s^{cc}_5(0.25),1) = 0.5625.
\end{split}\label{eqs:aff1}
\end{align}
With \eqref{eqs:aff1} we can improve the natural interval extensions range bounds $I_{5,nat}$ from $[-1.5,1]$ to $I_{5,alg} = [-0.75,0.5625]$ for factor $f_5$. The factor $f_5$ together with its convex and concave McCormick relaxations (constructed with natural interval bounds), natural interval bounds $f_5^L,f_5^U$ and the affine functions can be seen in Fig.~\ref{Fig:alg:interval(a)}. This procedure can be rerun for a different point in order to possibly improve the interval bounds even further but to keep this example simple, we do only one iteration. Note that it is possible to rerun the procedure for the same point if we improve a factor before, since the relaxations and the corresponding subgradients change.
\begin{figure}[t]
\centering
\begin{subfigure}{0.49\textwidth}
\begin{overpic}[width=\textwidth]{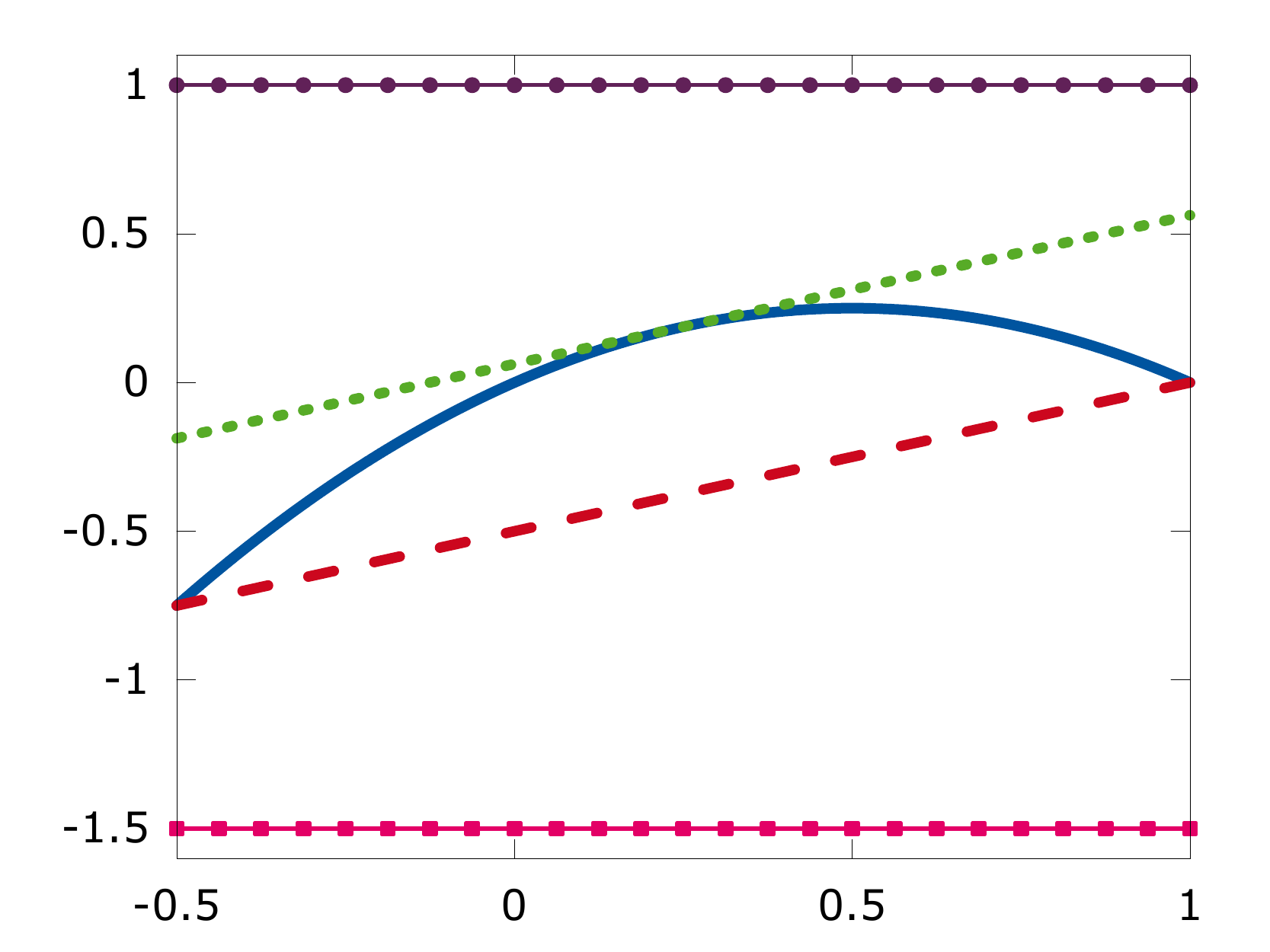}
\put (12.5,45) {\begin{tikzpicture}
\node[draw=none] at (0.36,0) {\tiny$f_5=f_5^{cc}$};
\node[draw=none] at (0.98,-0.3) {\tiny$f_5^{cv}=f_5^{cv,sub}(0.25)$};
\node[draw=none] at (0.61,-0.6) {\tiny$f_5^{cc,sub}(0.25)$};
\node[draw=none] at (0.04,-0.9) {\tiny $f_5^L$};
\node[draw=none] at (0.04,-1.2) {\tiny $f_5^U$};
\draw[line width=1pt,rwth, line cap = round] (-0.25,0)--(-0.7,0);
\draw[line width=1pt,dash pattern={on 5pt off 3pt },rot, line cap = round] (-0.25,-0.3)--(-0.7,-0.3);
\draw[line width=1pt,dash pattern={on 1pt off 3pt },grun, line cap = round] (-0.25,-0.6)--(-0.7,-0.6);
\draw[line width=1pt,magenta, line cap = round] (-0.25,-0.9)--plot[mark=square*,magenta,mark size=1.5pt] coordinates{(-0.475,-0.9)} -- (-0.7,-0.9);
\draw[line width=1pt,violett, line cap = round] (-0.25,-1.2) -- plot[mark=*,violett,mark size=1.5pt] coordinates{(-0.475,-1.2)} --(-0.7,-1.2);
\end{tikzpicture}}
\put (50,0) {\begin{tikzpicture}
\node[draw=none] at (0,0) {\small$z$};
\end{tikzpicture}}
\end{overpic}
\caption{}\label{Fig:alg:interval(a)}
\end{subfigure}%
\begin{subfigure}{0.49\textwidth}
\begin{overpic}[width=\textwidth]{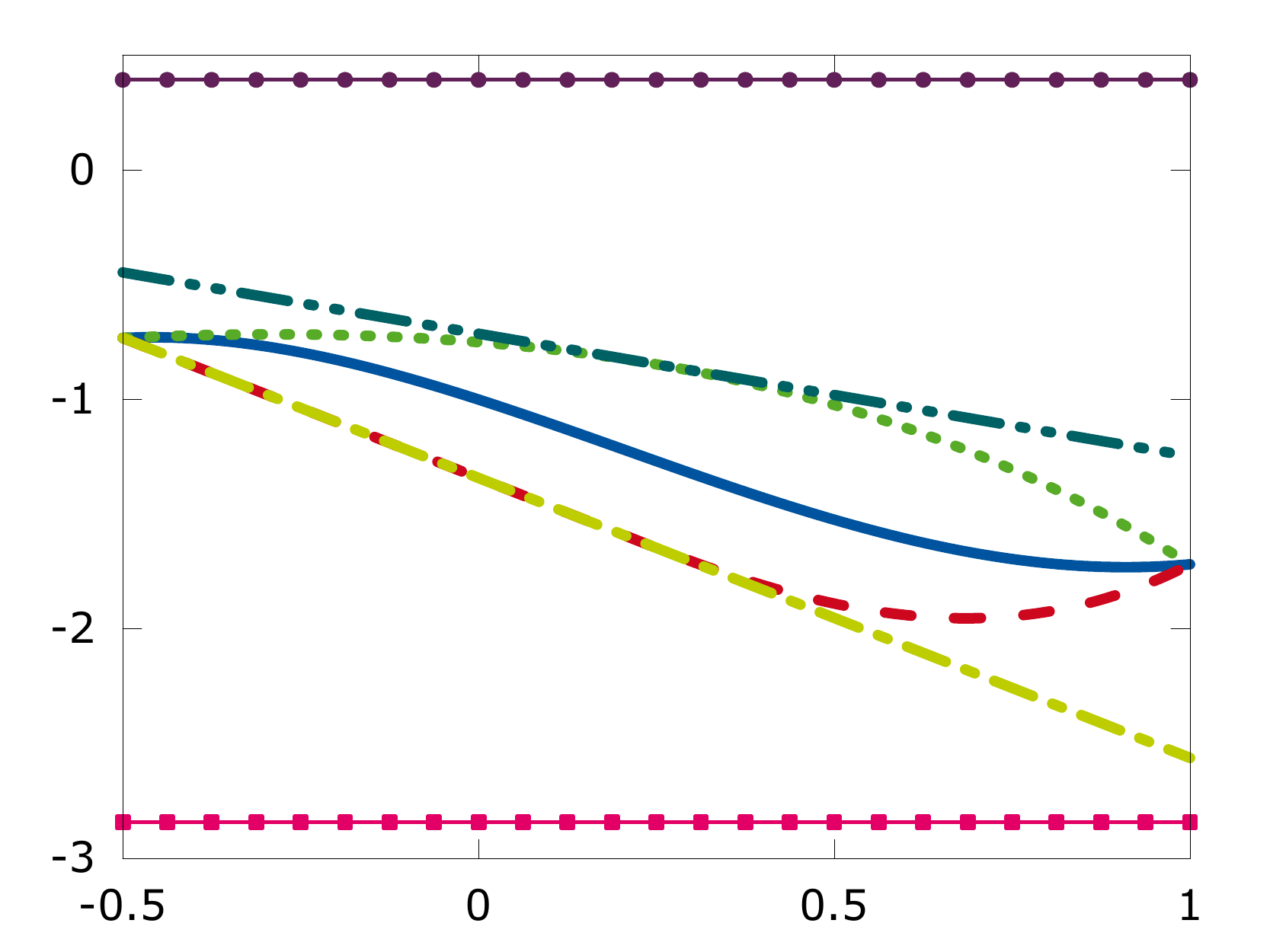}
\put (51,48) {\begin{tikzpicture}
\node[draw=none] at (0.3,0) {\tiny$f_6$};
\node[draw=none] at (0.24,-0.25) {\tiny$f_6^{cv}$};
\node[draw=none] at (0.24,-0.5) {\tiny$f_6^{cc}$};
\node[draw=none] at (-0.32,-0.75) {\tiny$f_6^{cv,sub}(0.25)$};
\node[draw=none] at (-0.32,-1) {\tiny$f_6^{cc,sub}(0.25)$};
\node[draw=none] at (-0.65,0) {\tiny $f_6^L$};
\node[draw=none] at (-0.65,-0.25) {\tiny $f_6^U$};
\draw[line width=1pt,rwth, line cap = round] (0.5,0)--(0.95,0);
\draw[line width=1pt,dash pattern={on 5pt off 3pt },rot, line cap = round] (0.5,-0.25)--(0.95,-0.25);
\draw[line width=1pt,dash pattern={on 1pt off 3pt },grun, line cap = round] (0.5,-0.5)--(0.95,-0.5);
\draw[line width=1pt,dash pattern={on 1pt off 3pt on 5pt off 3pt},maigrun, line cap = round] (0.5,-0.75)--(0.95,-0.75);
\draw[line width=1pt,dash pattern={on 5pt off 3pt on 1pt off 3pt on 1pt off 3pt},petrol, line cap = round] (0.5,-1)--(0.95,-1);
\draw[line width=1pt,magenta, line cap = round] (-0.45,0)--plot[mark=square*,magenta,mark size=1.5pt] coordinates{(-0.225,0)} -- (0,0);
\draw[line width=1pt,violett, line cap = round] (-0.45,-0.25) -- plot[mark=*,violett,mark size=1.5pt] coordinates{(-0.225,-0.25)} --(0,-0.25);
\end{tikzpicture}}
\put (50,0) {\begin{tikzpicture}
\node[draw=none] at (0,0) {\small$z$};
\end{tikzpicture}}
\end{overpic}
\caption{}\label{Fig:alg:interval(b)}
\end{subfigure}%
\caption{Example \ref{ex:heuristic}. \textbf{(a)} Factor $f_5(z)=z-z^2$ with its convex and concave McCormick relaxations, natural interval extension estimators $f_5^L,f_5^U$ for the range of $f_5$ on $Z=[-0.5,1]$ and affine under- and overestimators constructed with the use of subgradients at the middle point $0.25$. The affine underestimator equals the convex relaxation of $f_5$.\newline
\textbf{(b)} Factor $f_6(z)=z^3-\exp(z)$ with its convex and concave McCormick relaxations, natural interval extension estimators $f_6^L,f_6^U$ for the range of $f_6$ on $Z=[-0.5,1]$ and affine under- and overestimators constructed with the use of subgradients at the middle point $0.25$. }
\end{figure}

Next, we compute improved range bounds for the factor $f_6(z)=z^3-\exp(z)$. The convex and concave McCormick relaxations of $f_6$ on $Z$ (using the supplementary material of \cite{scott2011generalized}) are given as
\begin{align*}
f_6^{cv}(z)&=\begin{cases} 
-0.125-\exp(-0.5)+\left(0.1875-\frac{\exp(-0.5)-\exp(1)}{-0.5-1}\right)(z+0.5)&\text{, for } z\leq 0.25 \\
z^3	-\exp(-0.5)	-\frac{\exp(-0.5)-\exp(1)}{-0.5-1}(z+0.5)&\text{, else} \end{cases}\\
f_6^{cc}(z)&=0.25+0.75z-\exp(z).
\end{align*}
The natural interval extensions provide $I_{6,nat}=[-\exp(1)-0.125,1-\exp(-0.5)]\approx[-2.843,0.393]$, while the exact range is given as $I_{6,e}\approx[-1.73,-0.728]$. The convex and concave subgradients of $f_6^{cv}$ and $f_6^{cc}$, respectively, at the middle point $0.25$ of $Z$ are $s_6^{cv}(0.25)=0.1875-\frac{\exp(-0.5)-\exp(1)}{-0.5-1}$ and $s^{cc}_6(0.25)=0.75-\exp(0.25)$. We construct the corresponding affine functions 
\begin{align*}
f_6^{cv,sub}(0.25,z)&=f_6^{cv}(0.25)+s^{cv}_6(0.25)(z-0.25) \text{ and }\\
 f_6^{cc,sub}(0.25,z)&=f_6^{cc}(0.25)+s^{cc}_6(0.25)(z-0.25).
\end{align*}
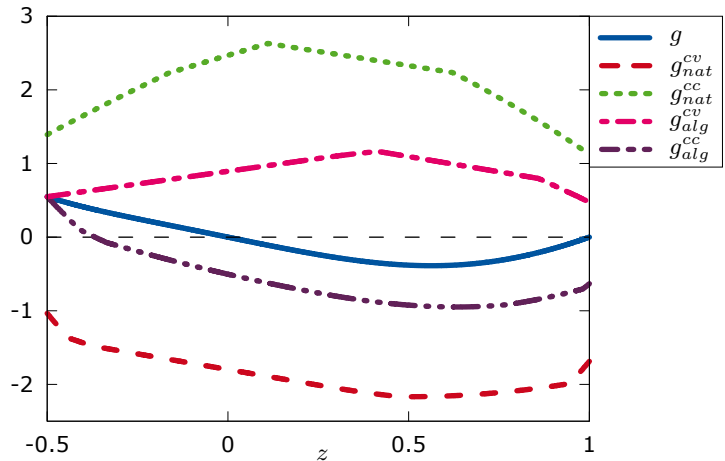
\begin{figure}[t]
\centering
\begin{overpic}[width=0.7\textwidth]{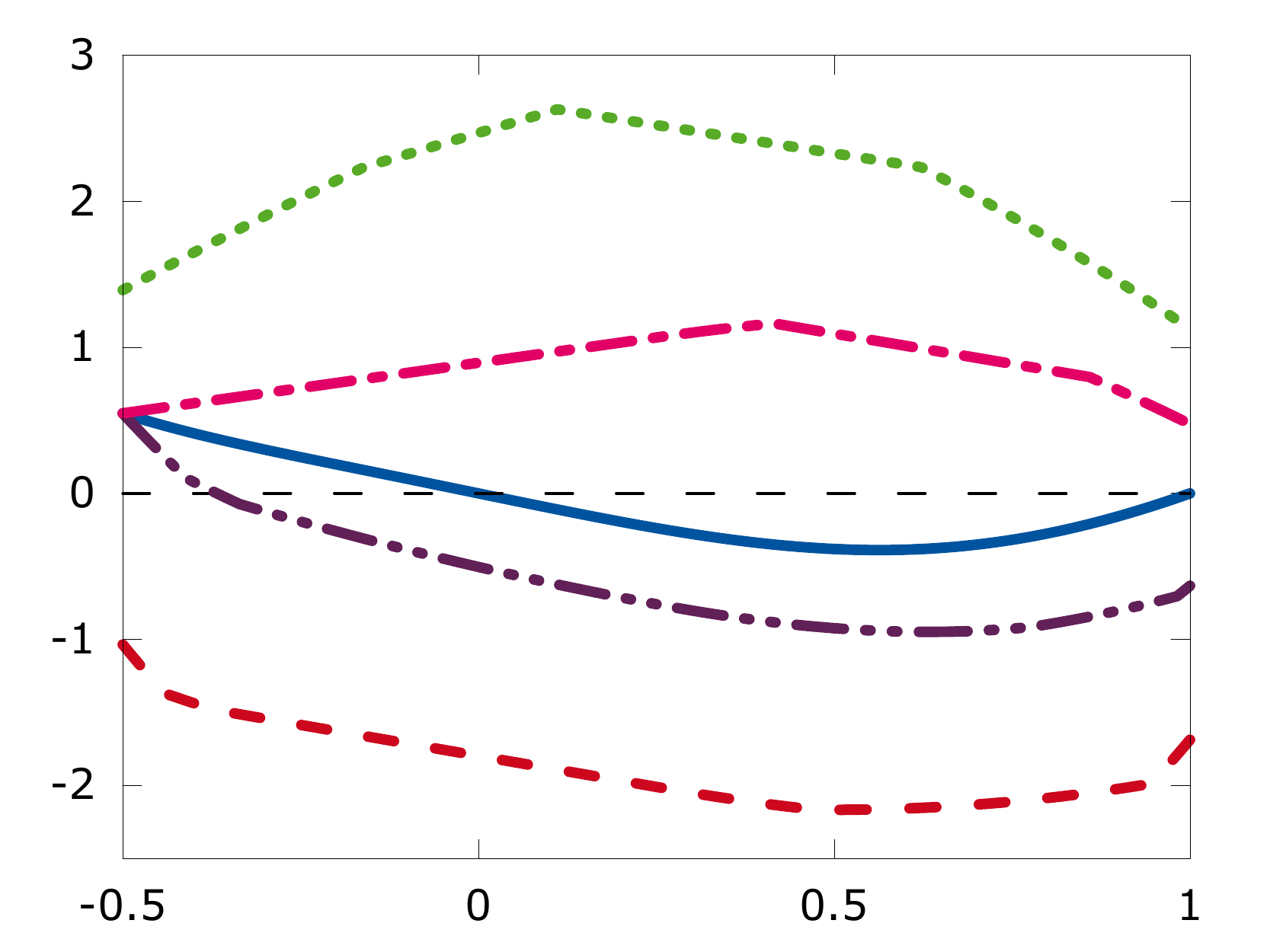}%
\put (95,46.5) {\begin{tikzpicture}
\node[draw=none] at (0,0) {\small$g$};
\node[draw=none] at (0.21,-0.37) {\small$g^{cv}_{nat}$};
\node[draw=none] at (0.21,-0.74) {\small$g^{cc}_{nat}$};
\node[draw=none] at (0.18,-1.11) {\small$g^{cv}_{alg}$};
\node[draw=none] at (0.18,-1.48) {\small$g^{cc}_{alg}$};
\draw[line width=2pt,rwth, line cap = round] (-0.35,0)--(-1,0);
\draw[line width=2pt,dash pattern={on 6pt off 6pt },rot, line cap = round] (-0.35,-0.37)--(-1,-0.37);
\draw[line width=2pt,dash pattern={on 1.5pt off 4pt },grun, line cap = round] (-0.35,-0.74)--(-1,-0.74);
\draw[line width=2pt,dash pattern={on 2pt off 4pt on 6pt off 4pt},magenta, line cap = round] (-0.35,-1.11)--(-1,-1.11);
\draw[line width=2pt,dash pattern={on 2pt off 4pt on 2pt off 4pt on 6pt off 4pt},violett, line cap = round] (-0.35,-1.48)--(-1,-1.48);
\end{tikzpicture}}
\put(93.7,47.04){\begin{tikzpicture}
\draw [line width=0.15pt] (0,0) rectangle (1.75,2); 
\end{tikzpicture}
}
\put (50,0) {\begin{tikzpicture}
\node[draw=none] at (0,0) {\small$z$};
\end{tikzpicture}}
\end{overpic}
\caption{Example \ref{ex:heuristic}. Function $g(z)=(z-z^2)(z^3-\exp(z))$ on $Z=[-0.5,1]$ together with its convex and concave McCormick relaxations $g^{cv}_{nat},g^{cc}_{nat}$ constructed with natural interval extensions and the convex and concave McCormick relaxation $g^{cv}_{alg},g^{cc}_{alg}$ constructed using the range bounds computed via heuristic presented in Algorithm \ref{alg:heuristic} with only 1 iteration at each factor.}
\label{Fig:alg:interval2}
\end{figure} 
Then, we evaluate the affine functions at their respective minimum and maximum in order to obtain
\begin{align}
\begin{split}
\min_{z\in [-0.5,1]} f_6^{cv,sub}(0.25,z) &= f_6^{cv,sub}(0.25,1) \approx -2.562\text{ and } \\
\max_{z\in [-0.5,1]} f_6^{cc,sub}(0.25,z) &= f_6^{cc,sub}(0.25,-0.5) \approx -0.446.
\end{split}\label{eqs:aff2}
\end{align}
With \eqref{eqs:aff2} we can improve the natural interval extensions range bounds $I_{6,nat}$ from $\approx[-2.843,0.393]$ to $I_{6,alg}\approx[-2.562,-0.446]$ for factor $f_6$. Factor $f_6$ together with its convex and concave McCormick relaxations (constructed with natural interval bounds), natural interval bounds $f_6^L,f_6^U$ and the affine functions can be seen in Fig.~\ref{Fig:alg:interval(b)}. Once again, this procedure can be rerun for a different point in order to possibly improve the interval bounds even further but to keep this example simple, we do only one iteration.

We can now construct the envelope for the bilinear product $f_7=f_5f_6$ on the improved intervals $I_{5,alg}\times I_{6,alg}$ and finally the convex and concave McCormick relaxations for $g$. Figure \ref{Fig:alg:interval2} shows function $g$ together with its convex and concave McCormick relaxations constructed with the simple intervals $I_{5,nat},I_{6,nat}$ denoted as $g^{cv}_{nat},g^{cc}_{nat}$ and two improved McCormick relaxations constructed with the new intervals $I_{5,alg},I_{6,alg}$ denoted as $g^{cv}_{alg},g^{cc}_{alg}$. The proposed heuristic drastically improves the relaxations. 
\end{example}

The procedure for obtaining improved lower and upper bounds for a factor of $g$ is formally given in Section \ref{sect:formaldefinition} and discussed afterwards.

\subsubsection{Other interval arithmetics}\label{ssec:otherIntervals}

Although natural interval extensions are very common in optimization for their simplicity, robustness and extremely low computational times, the usage of more sophisticated interval arithmetics is often advisable. If better interval extensions are used for the construction of McCormick relaxations, the resulting under- and overestimators may be a lot tighter than relaxations constructed with natural interval extensions, cf. Example 5 in \cite{Najman2017anchoring}. Still, even the better interval arithmetics do not guarantee that the resulting range bounds are exact. This leaves room for improvement of the range bounds by the presented heuristic. 

\begin{example}\label{ex:taylorform}
Consider $g(z)=(\log(z+1)-z^2)(\log(z+1)-\exp(z-0.5))$ on $Z=[-0.5,1]$, Fig.~\ref{Fig:alg2:interval2}, and consider the particular three factors $f_9(z)=\log(z+1)-z^2,f_{10}(z)=\log(z+1)-\exp(z-0.5)$ and $f_{11}(z)=f_9(z)\cdot f_{10}(z)$. For the other factors $z,1,0.5,z+1,z-0.5,z^2,\log(z+1),\exp(z-0.5)$, envelopes are known and simple interval arithmetics provide exact range bounds. In this example, we use the second order Taylor form interval extensions (Section 3.7 in \cite{nla.cat-vn1029192}) instead of natural interval extension to show that the presented algorithm can provide tighter McCormick relaxations even if more advanced interval arithmetics are used for the construction of McCormick relaxations. In particular, we compute the range bounds for a twice differentiable function $h:Z\to \mathbb{R}$ by calculating $I_{h,T}=[h_T^L,h_T^U]=h(c)+h'(c)(Z-c)+\frac{h''(Z)}{2}(Z-c)^2$, where $c$ is the middle point of $Z$, $h',h''$ are the first and second derivatives of $h$ and $h''(Z)$ is an interval overestimating the range of $h''$ which we calculated through natural interval extensions in this example.  

The McCormick relaxations of $f_9$ on $Z$ (using the supplementary material of \cite{scott2011generalized}) are given as
\begin{align*}
f_9^{cv}(z)&=f_9(-0.5)+\frac{f_9(-0.5)-f_9(1)}{-0.5-1}(z+0.5)\\
f_9^{cc}(z)&=\log(z+1)-z^2.
\end{align*}
The interval extensions obtained by the Taylor form provide $I_{9,T} \approx[-1.751,0.385]$ while the exact range is given as $I_{9,e}\approx[-0.943,0.177]$. The convex and concave subgradients of $f_9^{cv}$ and $f_9^{cc}$, respectively, at the middle point $0.25$ of $Z$ are $s_9^{cv}(0.25)=\frac{f_9^{cv}(-0.5)-f_9^{cv}(1)}{-0.5-1}$ and $s_9^{cc}(0.25)=0.3$, respectively. We construct the affine functions
\begin{align*}
f_9^{cv,sub}(0.25,z)&=f_9^{cv}(0.25)+s^{cv}_9(0.25)(z-0.25) \text{ and }\\
f_9^{cc,sub}(0.25,z)&=f_9^{cc}(0.25)+s^{cc}_9(0.25)(z-0.25).
\end{align*}
Subsequently, we compute the respective minimum and maximum of the affine functions
\begin{align}
\begin{split}
\min_{z\in [-0.5,1]} f_9^{cv,sub}(0.25,z) &= f_9^{cv,sub}(0.25,-0.5) \approx -0.943\text{ and } \\
\max_{z\in [-0.5,1]} f_9^{cc,sub}(0.25,z) &= f_9^{cc,sub}(0.25,1) \approx 0.385.
\end{split}\label{eqs:aff3}
\end{align}
With \eqref{eqs:aff3} we can improve the interval bounds obtained by the Taylor form from $I_{9,T}=[-1.751,0.385]$ to $I_{9,alg}=[-0.943,0.385]$ for factor $f_9$. Factor $f_9$ together with its convex and concave relaxations (constructed with Taylor form bounds), Taylor form bounds $f_{9,T}^L, f_{9,T}^U$ and the affine functions can be seen in Fig.~\ref{Fig:alg2:interval(a)}. Similar to Example \ref{ex:heuristic}, this procedure can be rerun but to keep this example simple, we do only one iteration.
\begin{figure}[t]
\centering
\begin{subfigure}{0.49\textwidth}
\begin{overpic}[width=\textwidth]{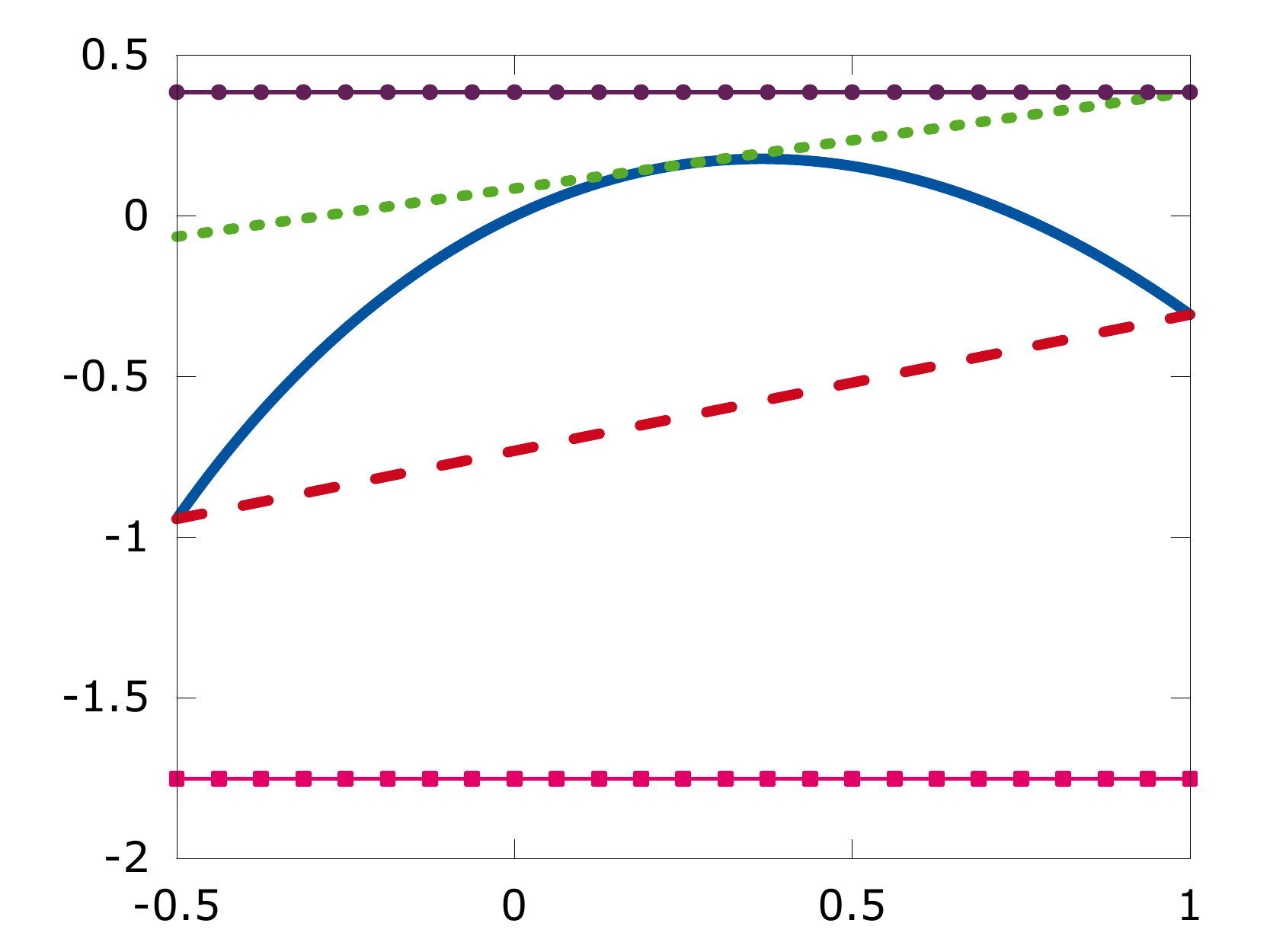}
\put (37,12) {\begin{tikzpicture}
\node[draw=none] at (-0.07,0) {\tiny$f_9=f_9^{cc}$};
\node[draw=none] at (-0.7,-0.3) {\tiny$f_9^{cv}=f_9^{cv,sub}(0.25)$};
\node[draw=none] at (-0.33,-0.6) {\tiny$f_9^{cc,sub}(0.25)$};
\node[draw=none] at (0.17,-0.9) {\tiny $f_{9,T}^L$};
\node[draw=none] at (0.17,-1.2) {\tiny $f_{9,T}^U$};
\draw[line width=1pt,rwth, line cap = round] (0.5,0)--(0.95,0);
\draw[line width=1pt,dash pattern={on 5pt off 3pt },rot, line cap = round] (0.5,-0.3)--(0.95,-0.3);
\draw[line width=1pt,dash pattern={on 1pt off 3pt },grun, line cap = round] (0.5,-0.6)--(0.95,-0.6);
\draw[line width=1pt,magenta, line cap = round] (0.5,-0.9)--plot[mark=square*,magenta,mark size=1.5pt] coordinates{(0.725,-0.9)} -- (0.95,-0.9);
\draw[line width=1pt,violett, line cap = round] (0.5,-1.2) -- plot[mark=*,violett,mark size=1.5pt] coordinates{(0.725,-1.2)} --(0.95,-1.2);
\end{tikzpicture}}
\put (50,0) {\begin{tikzpicture}
\node[draw=none] at (0,0) {\small$z$};
\end{tikzpicture}}
\end{overpic}
\caption{}\label{Fig:alg2:interval(a)}
\end{subfigure}%
\begin{subfigure}{0.49\textwidth}
\begin{overpic}[width=\textwidth]{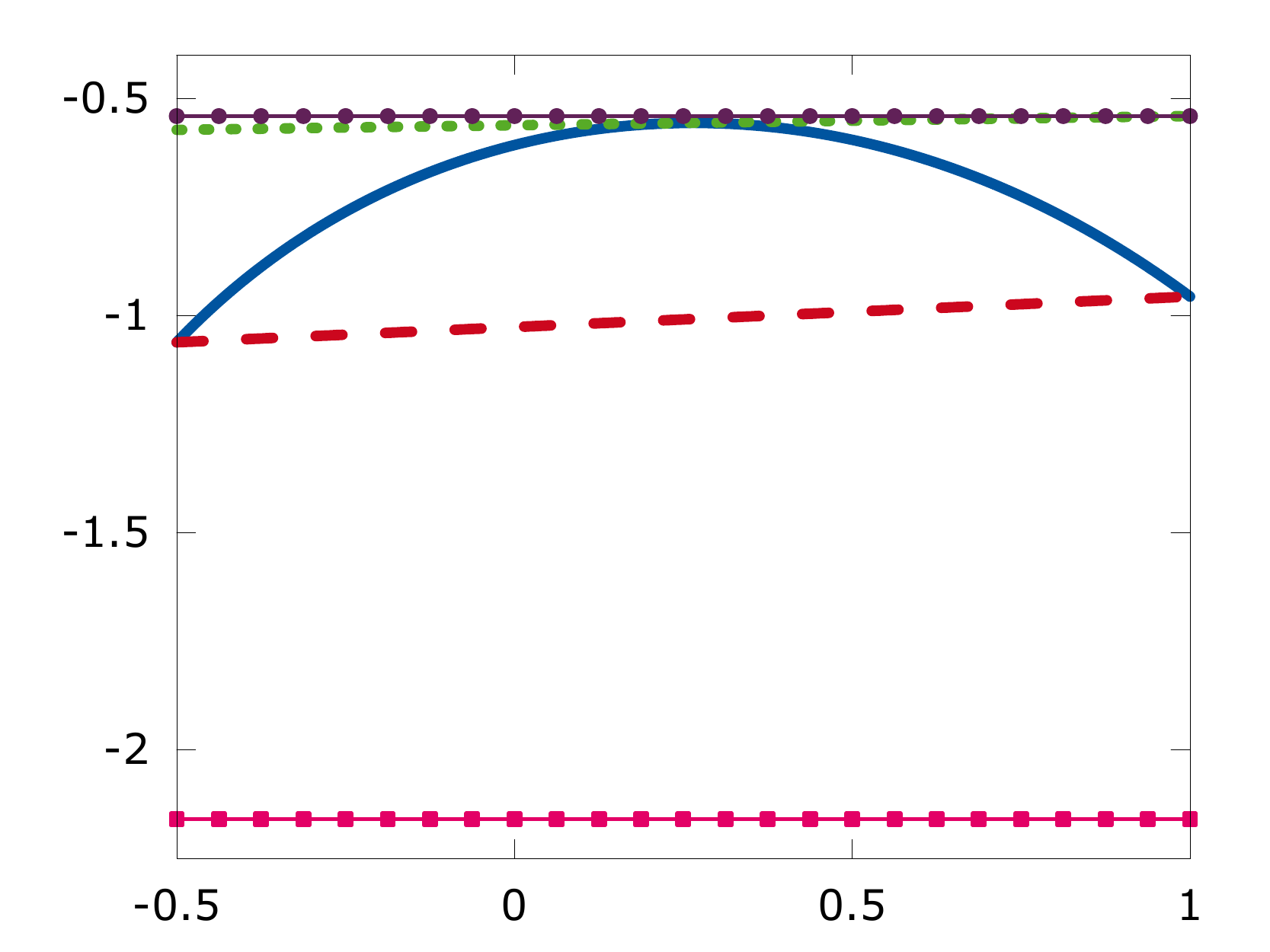}
\put (12.5,11) {\begin{tikzpicture}
\node[draw=none] at (0.42,0) {\tiny$f_{10}=f_{10}^{cc}$};
\node[draw=none] at (0.99,-0.3) {\tiny$f_{10}^{cv}=f_{10}^{cv,sub}(0.25)$};
\node[draw=none] at (0.61,-0.6) {\tiny$f_{10}^{cc,sub}(0.25)$};
\node[draw=none] at (0.19,-0.9) {\tiny $f_{{10},T}^L$};
\node[draw=none] at (0.20,-1.2) {\tiny $f_{{10},T}^U$};
\draw[line width=1pt,rwth, line cap = round] (-0.25,0)--(-0.7,0);
\draw[line width=1pt,dash pattern={on 5pt off 3pt },rot, line cap = round] (-0.25,-0.3)--(-0.7,-0.3);
\draw[line width=1pt,dash pattern={on 1pt off 3pt },grun, line cap = round] (-0.25,-0.6)--(-0.7,-0.6);
\draw[line width=1pt,magenta, line cap = round] (-0.25,-0.9)--plot[mark=square*,magenta,mark size=1.5pt] coordinates{(-0.475,-0.9)} -- (-0.7,-0.9);
\draw[line width=1pt,violett, line cap = round] (-0.25,-1.2) -- plot[mark=*,violett,mark size=1.5pt] coordinates{(-0.475,-1.2)} --(-0.7,-1.2);
\end{tikzpicture}}
\put (50,0) {\begin{tikzpicture}
\node[draw=none] at (0,0) {\small$z$};
\end{tikzpicture}}
\end{overpic}
\caption{}\label{Fig:alg2:interval(b)}
\end{subfigure}%
\caption{Example \ref{ex:taylorform}. \textbf{(a)} Factor $f_9(z)=\log(z+1)-z^2$ with its convex and concave McCormick relaxations, Taylor form interval extension estimators $f_{9,T}^L,f_{9,T}^U$ for the range of $f_9$ on $Z=[-0.5,1]$ and affine under- and overestimators constructed with the use of subgradients at the middle point $0.25$. The affine underestimator equals the convex relaxation of $f_9$.\newline
\textbf{(b)} Factor $f_{10}(z)=\log(z-1)-\exp(z-0.5)$ with its convex and concave McCormick relaxations, Taylor form interval extension estimators $f_{{10},T}^L,f_{{10},T}^U$ for the range of $f_{10}$ on $Z=[-0.5,1]$ and affine under- and overestimators constructed with the use of subgradients at the middle point $0.25$. The affine underestimator equals the convex relaxation of $f_{10}$. }
\end{figure}

Next, we calculate improved range bound for factor $f_{10}$. The convex and concave McCormick relaxations of $f_{10}$ on $Z$ (using the supplementary material of \cite{scott2011generalized}) are
\begin{align*}
f_{10}^{cv}(z)&=f_{10}(-0.5)+\frac{f_{10}(-0.5)-f_{10}(1)}{-0.5-1}(z+0.5) \\
f_{10}^{cc}(z)&=\log(z+1)-\exp(z-0.5).
\end{align*}
The Taylor form interval extensions provide $I_{{10},T}\approx[-2.16,-0.539]$, while the exact range is given as $I_{{10},e}\approx[-1.061,-0.555]$. The convex and concave subgradients of $f_{10}^{cv}$ and $f^{cc}_{10}$, respectively, at the middle point $0.25$ of $Z$ are $s_{10}^{cv}(0.25)=\frac{f_{10}(-0.5)-f_{10}(1)}{-0.5-1}$ and $s_{10}^{cc}(0.25)=0.8-\exp(-0.25)$. We construct the corresponding affine functions
\begin{align*}
f_{10}^{cv,sub}(0.25,z)&=f_{10}^{cv}(0.25)+s^{cv}_{10}(0.25)(z-0.25) \text{ and }\\
f_{10}^{cc,sub}(0.25,z)&=f_{10}^{cc}(0.25)+s^{cc}_{10}(0.25)(z-0.25)
\end{align*}
and compute the respective minimum and maximum
\begin{align}
\begin{split}
\min_{z\in [-0.5,1]} f_{10}^{cv,sub}(0.25,z) &= f_{10}^{cv,sub}(0.25,-0.5) \approx -1.061\text{ and } \\
\max_{z\in [-0.5,1]} f_{10}^{cc,sub}(0.25,z) &= f_{10}^{cc,sub}(0.25,1) \approx -0.539.
\end{split}\label{eqs:aff4}
\end{align}
With \eqref{eqs:aff4} we can improve the interval bounds obtained by the Taylor form $I_{{10},T}=\approx[-2.16,-0.539]$ to $I_{{10},alg}=[-1.061,-0.539]$ for factor $f_{10}$. Factor $f_{10}$ together with its convex and concave McCormick relaxations (constructed with Taylor form bounds), Taylor form bounds $f_{{10},T}^L,f_{{10},T}^U$ and the affine functions can be seen in Fig.~\ref{Fig:alg2:interval(b)}. Just as before, these steps can be recalculated but to keep this example simple, we do only one iteration.
\begin{figure}[t]
\centering
\begin{overpic}[width=0.7\textwidth]{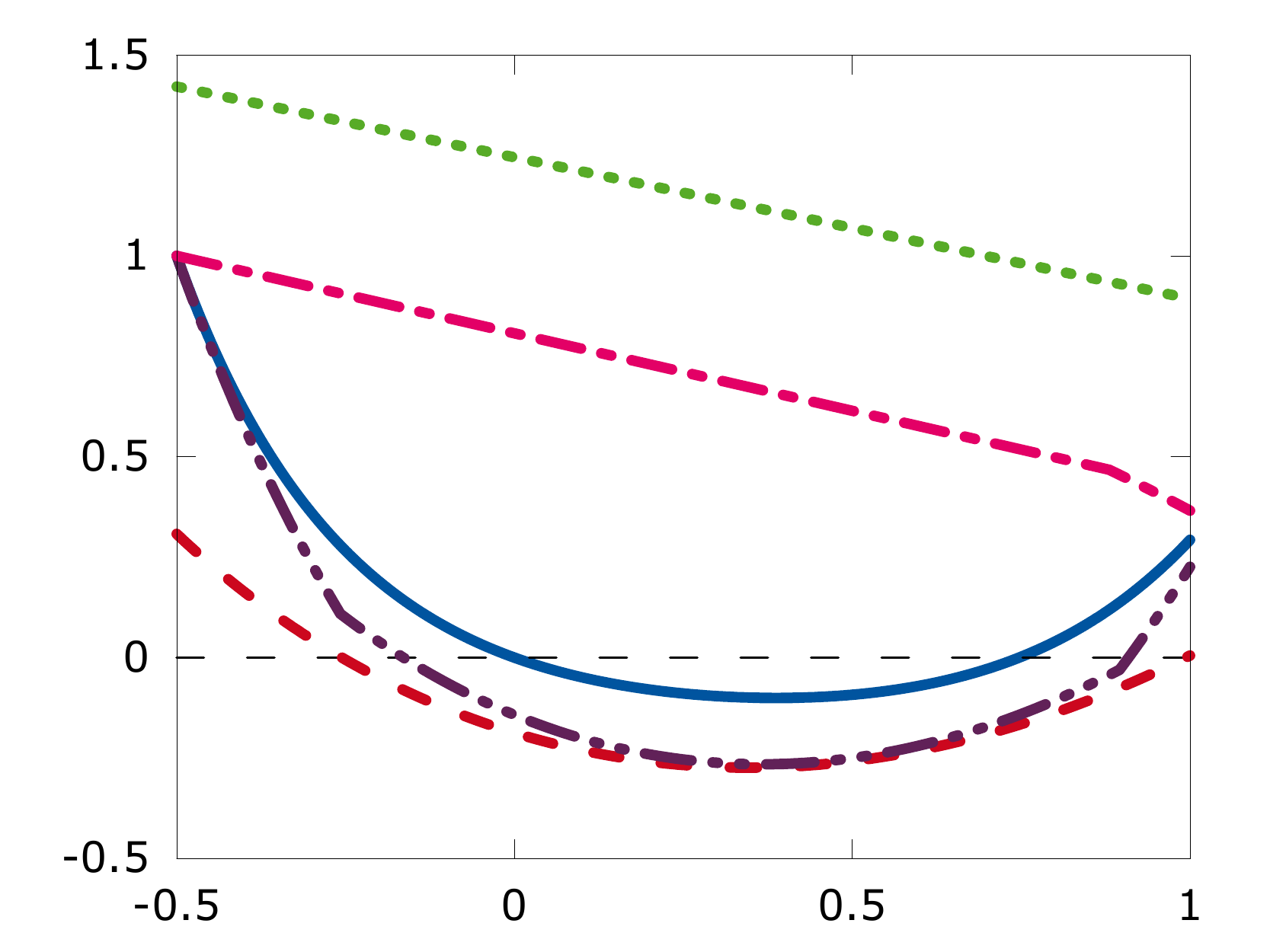}%
\put (95,46.5) {\begin{tikzpicture}
\node[draw=none] at (0,0) {\small$g$};
\node[draw=none] at (0.12,-0.37) {\small$g^{cv}_{T}$};
\node[draw=none] at (0.12,-0.74) {\small$g^{cc}_{T}$};
\node[draw=none] at (0.18,-1.11) {\small$g^{cv}_{alg}$};
\node[draw=none] at (0.18,-1.48) {\small$g^{cc}_{alg}$};
\draw[line width=2pt,rwth, line cap = round] (-0.35,0)--(-1,0);
\draw[line width=2pt,dash pattern={on 6pt off 6pt },rot, line cap = round] (-0.35,-0.37)--(-1,-0.37);
\draw[line width=2pt,dash pattern={on 1.5pt off 4pt },grun, line cap = round] (-0.35,-0.74)--(-1,-0.74);
\draw[line width=2pt,dash pattern={on 2pt off 4pt on 6pt off 4pt},magenta, line cap = round] (-0.35,-1.11)--(-1,-1.11);
\draw[line width=2pt,dash pattern={on 2pt off 4pt on 2pt off 4pt on 6pt off 4pt},violett, line cap = round] (-0.35,-1.48)--(-1,-1.48);
\end{tikzpicture}}
\put(93.7,47.04){\begin{tikzpicture}
\draw [line width=0.15pt] (0,0) rectangle (1.75,2); 
\end{tikzpicture}
}
\put (50,2) {\begin{tikzpicture}
\node[draw=none] at (0,0) {\small$z$};
\end{tikzpicture}}
\end{overpic}
\caption{Example \ref{ex:taylorform}. Function $g(z)=(\log(z+1)-z^2)(\log(z+1)-\exp(z-0.5))$ on $Z=[-0.5,1]$ together with its convex and concave McCormick relaxations $g^{cv}_{T},g^{cc}_{T}$ constructed with Taylor form interval extensions and the convex and concave McCormick relaxation $g^{cv}_{alg},g^{cc}_{alg}$ constructed using the range bounds computed via heuristic presented in Algorithm \ref{alg:heuristic} with only 1 iteration at each factor.}
\label{Fig:alg2:interval2}
\end{figure} 

Now, we are able to construct the envelope of $f_{11}=f_9f_{10}$ on $I_{9,alg}\times I_{{10},alg}$ and subsequently the McCormick relaxations of $g$ on $Z$. Figure~\ref{Fig:alg2:interval2} shows function $g$ together with its convex and concave McCormick relaxations constructed with the intervals $I_{9,T},I_{{10},T}$ obtained with the Taylor form interval extensions denoted as $g_{T}^{cv},g_T^{cc}$ and two improved McCormick relaxations constructed with the intervals $I_{9,alg},I_{{10},alg}$ denoted as $g_{alg}^{cv},g_{alg}^{cc}$. We see that even for the more sophisticated interval extensions, the proposed heuristic is able to significantly improve the final resulting McCormick relaxations. 
\end{example} 

Example \ref{ex:taylorform} shows that the heuristic is able to improve the McCormick relaxations of a given function $g$ even if more advanced interval arithmetics are used for the computation of range bounds of each factor. Obviously, it holds that the weaker the underlying estimated bounds for each factor, the larger is the potential of the presented heuristic.

\subsection{Formal definition of the algorithm}\label{sect:formaldefinition}

For a given factorable function $g:Z\to\mathbb{R},Z\in\mathbb{I}\mathbb{R}^n$, we traverse the corresponding DAG of $g$ starting at the independent variables $z_i,i\in\{1,\dots,n\}$ and working through all factors up to the root given as $g$ (reverse-level-order in Graph Theory terminology). In each factor $f_j,j \in \{1,\dots, |\mathbb{F}| \}$, we execute Algorithm \ref{alg:heuristic} in order to obtain bounds $I_{j,alg}=[f_{j,alg}^L,f_{j,alg}^U]$ on the range of $f_j$ and save these. We use $I_{j,alg}$ then directly when computing $f_{k},k>j$. Note that the computed range bounds for every $f_j$ are valid on whole $Z$. Thus, after the DAG of $g$ has been completely traversed, we can calculate McCormick relaxations and its subgradients of $g$ at any point $\bar{\vec{z}}\in Z$ with the use of the range bounds $I_{j,alg}$ for each factor $f_j,j\in\{1,\dots,|\mathbb{F}|\}$ instead of using natural interval extensions for the range bounds estimation.  

\begin{algorithm}
Given a DAG representation $G=(\mathbb{F},E)$ of a factorable function $g:Z\to\mathbb{R},~Z\in\mathbb{I}\mathbb{R}^n$. \newline
$\vec{z}^L,\vec{z}^U$ - lower and upper, finite, bounds for all independent variables $z_i, ~i\in \{1,\dots,n \}$\newline
$\bar{\vec{z}}$ - initial point for computation of relaxations and subgradients \newline
$N$ - maximal number of iterations for given factor\newline
Initialize factor $f_j \in \mathbb{F}$ obtained by traversing $G$ in reversed-level-order to obtain $f_j^{cv}(\bar{\vec{z}}),f_j^{cc}(\bar{\vec{z}}),\vec{s}_j^{cv}(\bar{\vec{z}}),\vec{s}_j^{cc}(\bar{\vec{z}}),I_j=[f_j^L,f_j^U]$.  

k=1\;
\While{$k\leq N$}{ \label{line:while}
 \If{$f_j^L<f_j^U$}{ \label{line:check_const}
 
	 t$^{cv}$ = $f_j^{cv}(\bar{\vec{z}})$\; \label{line:start}
	 t$^{cc}$ = $f_j^{cc}(\bar{\vec{z}})$\;
	 \For{i=1,\dots,n}{ \label{line:for}
			 \uIf{$s_i^{cv}(\bar{\vec{z}})\geq 0$}{\label{line:eps}
				 t$^{cv}$=t$^{cv}$+$s_i^{cv}(\bar{\vec{z}})(z^L_i-\bar{z}_i)$\;
			 }\Else{
				 t$^{cv}$=t$^{cv}$+$s_i^{cv}(\bar{\vec{z}})(z^U_i-\bar{z}_i)$\;
			 }
			\uIf{$s_i^{cc}(\bar{\vec{z}})\geq 0$}{\label{line:eps2}
				 t$^{cc}$=t$^{cc}$+$s_i^{cc}(\bar{\vec{z}})(z^U_i-\bar{z}_i)$\;
			 }\Else{
				 t$^{cc}$=t$^{cc}$+$s_i^{cc}(\bar{\vec{z}})(z^L_i-\bar{z}_i)$\;
			 }
	 }
		$f_{j,alg}^L=\max\{$t$^{cv}, f_j^L \}$\;\label{line:if1}
		$f_{j,alg}^U=\min\{$t$^{cc}, f_j^U \}$\;\label{line:if2}
		$I_{j}=[f_{j,alg}^L,f_{j,alg}^U]$\; \label{line:newinterval}		
		Save $I_{j,alg}=I_{j}$\;	\label{line:save}
		\If{$k+1\leq N$}{ 
			Compute $\bar{\vec{z}}_{new}$\; \label{line:newpoint}
		}	
	}
 k=k+1\;
} 

\caption{Method for obtaining tighter interval range bounds with the use of propagated subgradients for one factor $f_j$ of $g$. }\label{alg:heuristic}
\end{algorithm}  

\subsection{Algorithm discussion}\label{subsec:disc}
We now discuss Algorithm \ref{alg:heuristic}. First, we check if the factor $f_j$, we currently consider, is a constant function by simply comparing the lower and upper bound (line \ref{line:check_const}), since if it is the case, we cannot improve any bounds on its range and thus, the heuristic is unnecessary. Obviously this check is only sufficient and not necessary, since it is possible to have different bounds for a constant function due to overestimation. In the \texttt{for} loop (line \ref{line:for}), we solve $\min\limits_{\vec{z}\in Z} f_j^{cv,sub}(\bar{\vec{z}},\vec{z})$ and $\max\limits_{\vec{z}\in Z} f_j^{cc,sub}(\bar{\vec{z}},\vec{z})$ by simple subgradient comparisons as both problems are box-constrained and linear. The correct corner $\vec{z}^c$ of $Z$ is determined by examination of the sign of the subgradient in the particular dimension (lines \ref{line:eps} and \ref{line:eps2}). Then, we check if we can improve the bounds of factor $f_j$ (lines \ref{line:if1} and \ref{line:if2} ). We update the range bounds $I_j$ of factor $f_j$ (line \ref{line:newinterval}) such that they can be directly used for the computation of relaxations, subgradients and bounds for factors $f_{k},k>j$. Then, we save the improved bounds (line \ref{line:save}). 
If needed, we compute a new point $\bar{\vec{z}}_{new}$ (line \ref{line:newpoint}). There are many ways to compute a new point. In our computational studies in Section \ref{sect:numeric}, we use the simple bisection in order to obtain a new $\bar{\vec{z}}_{new}$, i.e., the next point is given by the middle point of the interval $[\bar{\vec{z}}, \vec{z}^c]$, where $\vec{z}^c$ is computed before in lines \ref{line:eps} and \ref{line:eps2} (not explicitly shown in Algorithm \ref{alg:heuristic}). Note that the simple bisection method converges to $\min\limits_{\vec{z}\in Z} f_j^{cv}(\vec{z})$ for $N\to\infty$. A more sophisticated method for the computation of $\bar{\vec{z}}_{new}$ may provide better results and represents potential future work. Note that it also makes sense to compute $2$ new points in order to independently improve the upper and lower bound within the algorithm. 

We are also interested in the computational complexity of the presented heuristic. The computation of the possibly improved bounds $I_{j,alg}=[f_{j,alg}^L,f_{j,alg}^U]$ (lines \ref{line:start} - \ref{line:save}) is linear in the dimension of the optimization variables $\vec{z}$, i.e., we have a complexity of $\mathcal{O}(n)$ for the computation of $I_{j,alg}$. This is comparable to the computation of McCormick relaxations and propagation of subgradients which have to be computed for each dimension in each factor. If only 1 ($N=1$) iteration is allowed, we can directly use the propagated relaxation and subgradient values at the desired point $\bar{\vec{z}}$ and are not forced to re-evaluate McCormick relaxations, subgradients and interval bounds for all the factors on which $f_j$ depends for subsequent points. Therefore, if a function $g$ consists of $|\mathbb{F}|$ factors and we allow only one iteration of the heuristic, the computational complexity amounts to $\mathcal{O}(|\mathbb{F}|\cdot n)$. Regarding the computational time for $N=1$, we can expect that even in cases where the heuristic does not yield any improvement, the additional computational effort is negligible. This observation is also confirmed in the numerical studies in Section \ref{sect:numeric}. %
In contrast, if we allow for more than $1$ iteration ($N>1$) within a factor $f_j$, we have to propagate all required information through all previous factors at the new point. Indeed, if a factor $f_j$ depends on $\mathcal{J}$ other factors and we allow $N$ iterations, we need to do $\mathcal{J}^N$ computations, which is large for more complicated functions ($\mathcal{J}\gg 1$) and $N>1$. Moreover, let us assume that each factor depends on all the previous factors, then the complexity for $N$ computations of $f_j$ equals $\sum\limits_{k=1}^{j} k^N$, which is a polynomial complexity but still extremely large if a function consists of many factors and $N>1$. The impact of the number of factors for $N>1$ matches the numerical results presented in Section \ref{sect:numeric}. To avoid a large number of re-computations, we could heuristically decide whether it is worth to traverse the DAG again, e.g., by the value of $j$ or if the difference between the McCormick relaxation $f_j^{cv}$ and the natural interval bounds $f_j^L$ is very large.  

\subsection{Algorithm limitations}
The algorithm cannot deteriorate the bounds, since new bounds $I_{j,alg}$ are given by the maximum and minimum of the originally computed bounds $f_j^L,f_j^U$ and $t^{cv},t^{cc}$ (lines \ref{line:if1} and \ref{line:if2}). The best bounds obtained by Algorithm \ref{alg:heuristic} cannot be better than the minimum and maximum of the convex and concave McCormick relaxations of a factor $f_j$, i.e., it holds that
\begin{align*}
f^L_{j,alg} \leq \min_{\vec{z}\in Z} f^{cv}(\vec{z}) ~\text{ and }~ f^U_{j,alg} \geq \max_{\vec{z}\in Z} f^{cc}(\vec{z}).
\end{align*}  
Algorithm \ref{alg:heuristic} is able to improve the range bounds of a factor, see Example \ref{ex:heuristic} but is not guaranteed to improve the bounds of a factor $f_j$, e.g., if the underlying interval extensions for the bounds $f_j^L,f_j^U$ are already exact or if the point $\bar{\vec{z}}$ is chosen badly as we show in the next example.

\begin{example}\label{ex:heuristic2}
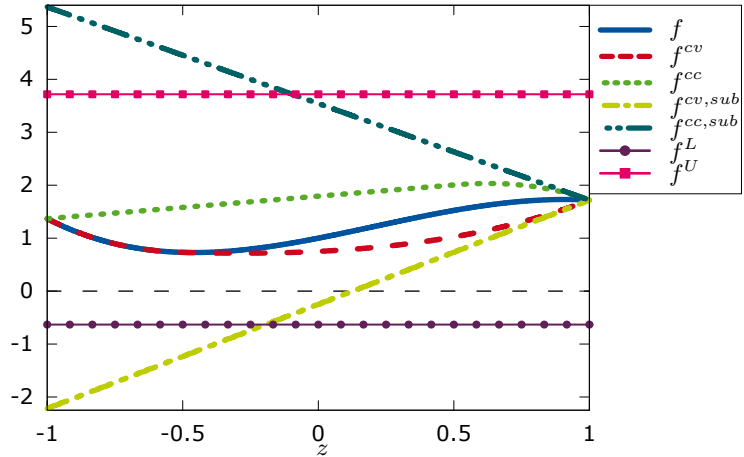
\begin{figure}[t]
\centering
\begin{overpic}[width=0.7\textwidth]{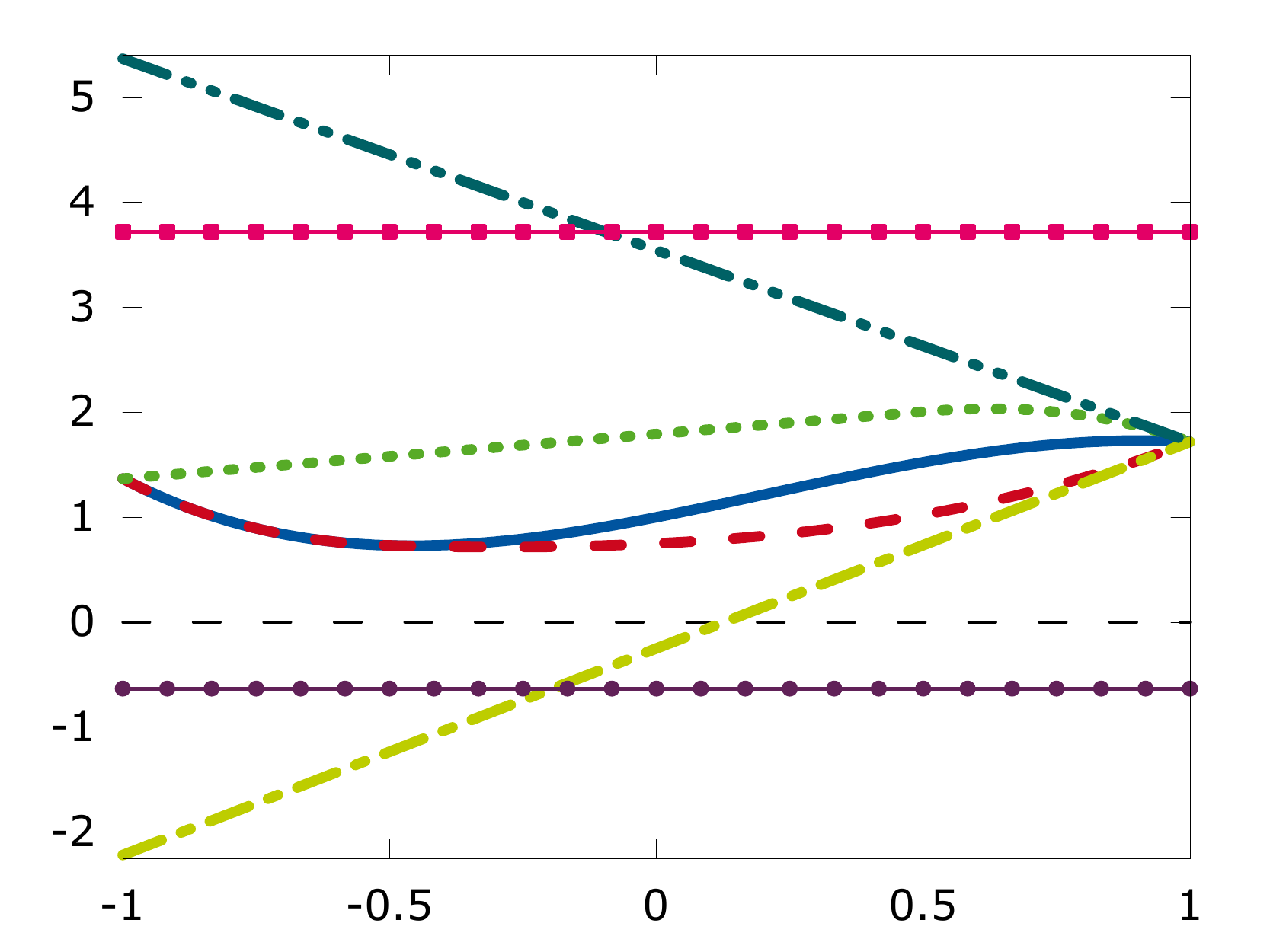}%
\put (95,41) {\begin{tikzpicture}
\node[draw=none] at (0,0) {\small$f$};
\node[draw=none] at (0.14,-0.35) {\small$f^{cv}$};
\node[draw=none] at (0.13,-0.7) {\small$f^{cc}$};
\node[draw=none] at (0.38,-1) {\small$f^{cv,sub}$};
\node[draw=none] at (0.37,-1.3) {\small$f^{cc,sub}$};
\node[draw=none] at (0.1,-1.6) {\small$f^L$};
\node[draw=none] at (0.11,-1.9) {\small$f^U$};
\draw[line width=2pt,rwth, line cap = round] (-0.35,0)--(-1,0);
\draw[line width=2pt,dash pattern={on 4pt off 3pt },rot, line cap = round] (-0.35,-0.35)--(-1,-0.35);
\draw[line width=2pt,dash pattern={on 1pt off 3pt },grun, line cap = round] (-0.35,-0.7)--(-1,-0.7);
\draw[line width=2pt,dash pattern={on 5pt off 3pt on 1pt off 3pt},maigrun, line cap = round] (-0.35,-1)--(-1,-1);
\draw[line width=2pt,dash pattern={on 8pt off 3pt on 1pt off 3pt on 1pt off 3pt},petrol, line cap = round] (-0.35,-1.3)--(-1,-1.3);
\draw[line width=1pt,magenta, line cap = round] (-0.35,-1.9)--plot[mark=square*,magenta,mark size=1.5pt] coordinates{(-0.675,-1.9)} -- (-1,-1.9);
\draw[line width=1pt,violett, line cap = round] (-0.35,-1.6) -- plot[mark=*,violett,mark size=1.5pt] coordinates{(-0.675,-1.6)} --(-1,-1.6);
\end{tikzpicture}}
\put(93.7,41.15){\begin{tikzpicture}
\draw [line width=0.15pt] (0,0) rectangle (2,2.5); 
\end{tikzpicture}
}
\put (50,-1) {\begin{tikzpicture}
\node[draw=none] at (0,0) {\small$z$};
\end{tikzpicture}}
\end{overpic}
\caption{Example \ref{ex:heuristic2}. Factor $f(z)=\exp(z)-z^3$ on $Z=[-1,1]$. The heuristic does not provide any improvement if the initial point is chosen as $\bar{z}=1$ and only $1$ iteration ($N=1$) is allowed. }
\label{Fig:alg:interval3}
\end{figure} 
Consider $f(z)=\exp(z)-z^3$ on $Z=[-1,1]$. Let us apply Algorithm \ref{alg:heuristic} to $f$ at point $\bar{z}=1$. The heuristic does not improve the interval bounds of $f$, see Fig.~\ref{Fig:alg:interval3}. 
\end{example}

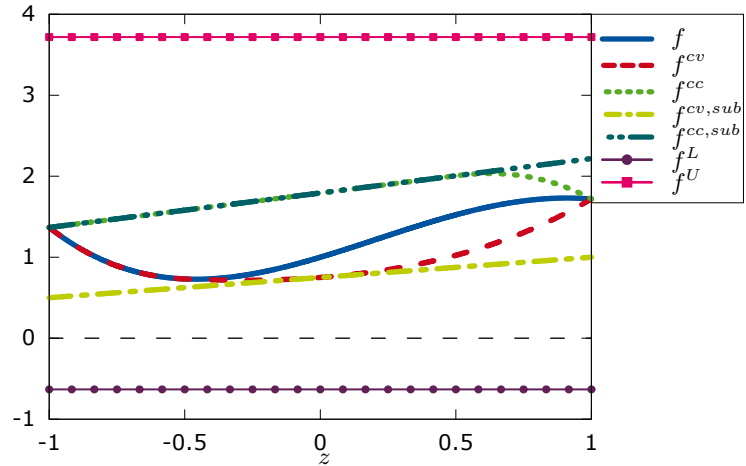
\begin{figure}[t]
\centering
\begin{overpic}[width=0.7\textwidth]{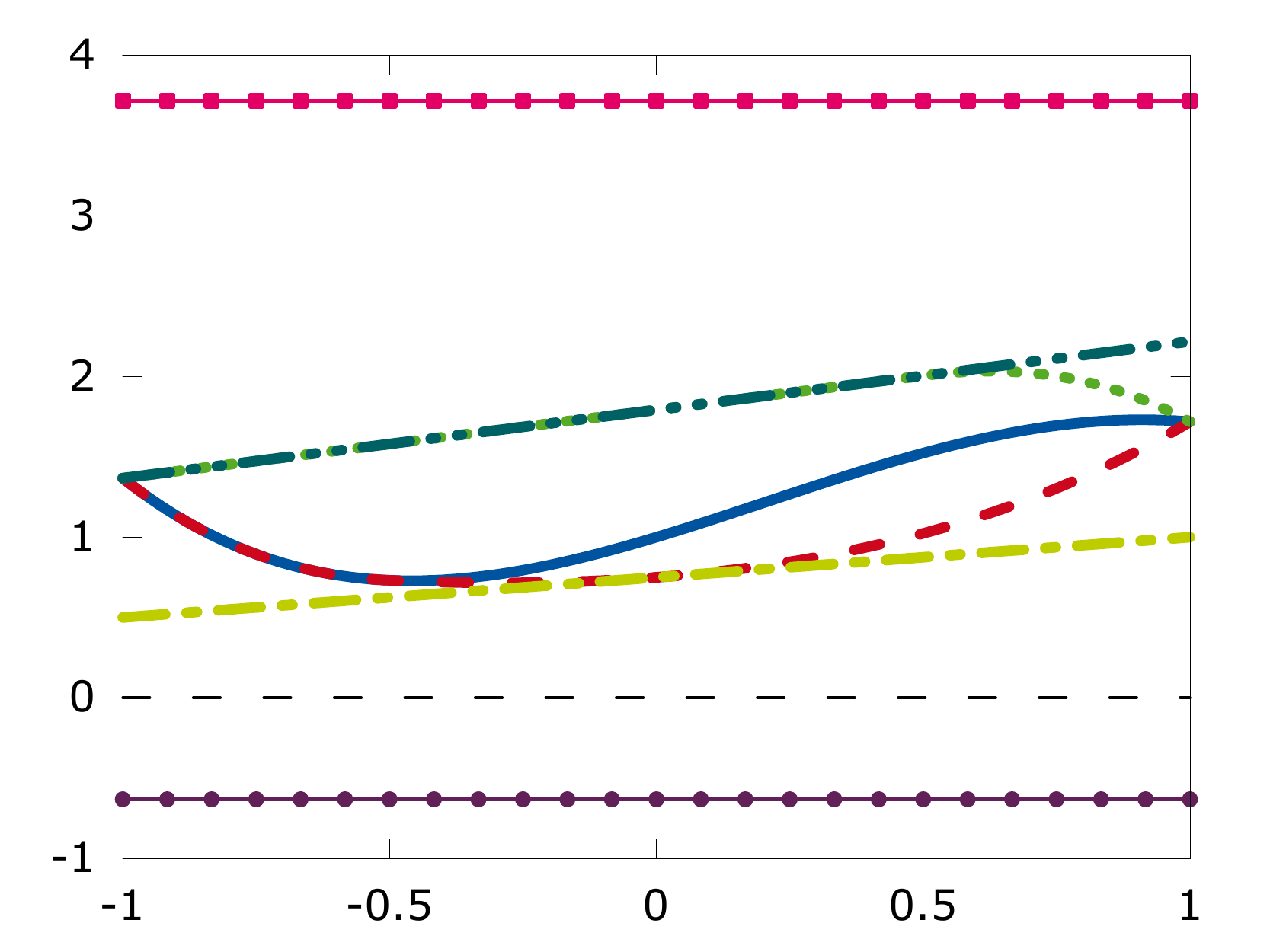}%
\put (95,41) {\begin{tikzpicture}
\node[draw=none] at (0,0) {\small$f$};
\node[draw=none] at (0.14,-0.35) {\small$f^{cv}$};
\node[draw=none] at (0.13,-0.7) {\small$f^{cc}$};
\node[draw=none] at (0.38,-1) {\small$f^{cv,sub}$};
\node[draw=none] at (0.37,-1.3) {\small$f^{cc,sub}$};
\node[draw=none] at (0.1,-1.6) {\small$f^L$};
\node[draw=none] at (0.11,-1.9) {\small$f^U$};
\draw[line width=2pt,rwth, line cap = round] (-0.35,0)--(-1,0);
\draw[line width=2pt,dash pattern={on 4pt off 3pt },rot, line cap = round] (-0.35,-0.35)--(-1,-0.35);
\draw[line width=2pt,dash pattern={on 1pt off 3pt },grun, line cap = round] (-0.35,-0.7)--(-1,-0.7);
\draw[line width=2pt,dash pattern={on 5pt off 3pt on 1pt off 3pt},maigrun, line cap = round] (-0.35,-1)--(-1,-1);
\draw[line width=2pt,dash pattern={on 8pt off 3pt on 1pt off 3pt on 1pt off 3pt},petrol, line cap = round] (-0.35,-1.3)--(-1,-1.3);
\draw[line width=1pt,magenta, line cap = round] (-0.35,-1.9)--plot[mark=square*,magenta,mark size=1.5pt] coordinates{(-0.675,-1.9)} -- (-1,-1.9);
\draw[line width=1pt,violett, line cap = round] (-0.35,-1.6) -- plot[mark=*,violett,mark size=1.5pt] coordinates{(-0.675,-1.6)} --(-1,-1.6);
\end{tikzpicture}}
\put(93.7,41.15){\begin{tikzpicture}
\draw [line width=0.15pt] (0,0) rectangle (2,2.5); 
\end{tikzpicture}
}
\put (50,-1) {\begin{tikzpicture}
\node[draw=none] at (0,0) {\small$z$};
\end{tikzpicture}}
\end{overpic}
\caption{Factor $f(z)=\exp(z)-z^3$ on $Z=[-1,1]$. The heuristic provides clear improvement if an additional iteration is allowed or if the initial point is directly set to $\bar{z}=0$ in Example \ref{ex:heuristic2}.}
\label{Fig:alg:interval4}
\end{figure}

Example \ref{ex:heuristic2} shows that the outcome of the heuristic depends on the chosen initial point $\bar{\vec{z}}$ and also on the maximum number of iterations. Choosing a corner point $\vec{z}^c\in Z$ as initial point for Algorithm \ref{alg:heuristic} is only a good choice if we have some monotonicity information of the convex relaxation. In general, choosing the initial point for the heuristic from the interior of $Z$ seems more intuitive and promising due to the positive curvature of convex underestimators. Note that the bounds of $f$ in Example \ref{ex:heuristic2} are improved if additional iterations of the heuristic are performed. If we allow for an additional iteration in Example \ref{ex:heuristic2}, we obtain the middle point $\bar{\vec{z}}_{new}=0$ by applying simple bisection of $Z$ for the re-computation of $\bar{\vec{z}}$ and the heuristic indeed does improve the range bounds of $f$, which can be seen in Fig.~\ref{Fig:alg:interval4}.

\section{Numerical results}\label{sect:numeric}

In order to test the presented heuristic, we use the software implementation given in Section 4 of \cite{bongartz2017deterministic}. In particular we use CPLEX v12.5 \cite{Cplex2009} for linear optimization, the IPOPT solver \cite{Waechter2006} for local nonlinear optimization, the FADBAD\text{++} package for automatic differentiation \cite{fadbad2016} and the MC\texttt{++} package v2.0 \cite{Chachuatmc++} for McCormick relaxations. All calculations were conducted on an Intel$^{\text{\textregistered}}$ Core\texttrademark i3-3240 with 3.4GHz and 8 GB RAM on Windows 7.  We solved 13 small problems of varying sizes with up to 14 variables chosen arbitrarily from the COCONUT library \cite{Shcherbina2003}; the well-known six-hump-camel function and the ursem-waves function; 3 case studies of a combined-cycle power plant presented in Section 5 of \cite{bongartz2017deterministic} for minimizing the levelized cost of electricity and the 2 case studies for parameter estimation presented in Section 5 of \cite{mitsos2009mccormick}. We solved the nonlinear optimization problems as explained in the following. The upper bounding problems were solved locally with the local solver IPOPT in order to obtain a valid upper bound. For the lower bound, we relaxed the problems using McCormick relaxations and then constructed linearizations $g^{cv,sub}(\bar{\vec{z}},\vec{z})$ at a single point $\bar{\vec{z}}\in Z$ of the convex McCormick relaxations with the use of subgradient propagation. The considered problems consist of constrained and box-constrained problems. In both cases we linearize the convex relaxation of the objective function and all constraints (excluding the variable bounds) at a predefined single point to construct a linear program, which is then solved with CPLEX. In the case of  box-constraints only, we obtain an extremely simple linear program of which the solution can be obtained by simple coefficient analysis. Still, even in this simple case we automatically call CPLEX for the solution. We always use only one linearization point, namely the middle point of the current node in the first and third numerical comparisons (Tables \ref{table2} and \ref{tableIterations}) and the current incumbent in the second comparison (Table \ref{table3}). If any of the coordinates of the incumbent is not within a given node, we simply replace it by the corresponding middle point coordinate. 

Additionally, in most problems, there were no bounds given for the optimization variables. Since McCormick relaxations need valid bounds in order to be constructed, we provided valid bounds containing the global minimum in all problems, see Table \ref{table1} in Appendix \ref{appendix} for the number of variables, inequalities, equalities and bounds we used. If the bounds are provided in the given optimization problem, we marked it with the keyword \emph{given}. Moreover, problem \texttt{ex6\_2\_10} consists of many uses of the convex function $f(x)=x\cdot \log(x)$. In this particular problem we used the \texttt{xlog} function, implemented in MC\texttt{++}(v2.0)\cite{Chachuatmc++} instead of using the binary product of $x$ and $\log(x)$. We also used the envelope of the logarithmic mean temperature difference function, presented in \cite{Najman20161605} for the 3 case studies from \cite{bongartz2017deterministic} explaining the improved computational times for $N=0$ in this article in comparison to \cite{bongartz2017deterministic}. 

In the following, we discuss the numerical impact of the heuristic within a B\&B framework using McCormick relaxations. We do not compare with other relaxation techniques such as the auxiliary variable method or $\alpha$BB \cite{androulakis_95_1,maranas1992global}, since we are interested in the improvement of McCormick relaxations in particular. Also, applying more sophisticated lower bounding procedures, e.g., using more than only one linearization point or applying a bundle method, may be of interest and represents potential future work but is out of the scope of this article. It is worth mentioning that the great potential of the described B\&B procedure in the sense of computational time compared to state-of-the-art deterministic global optimization solvers has already been shown in \cite{bongartz2017deterministic,mitsos2009mccormick} for numerical experiments \texttt{Case study I,II,III} \cite{bongartz2017deterministic}, \texttt{heat} and \texttt{kinetic} \cite{mitsos2009mccormick}. 

\begin{figure}[t]
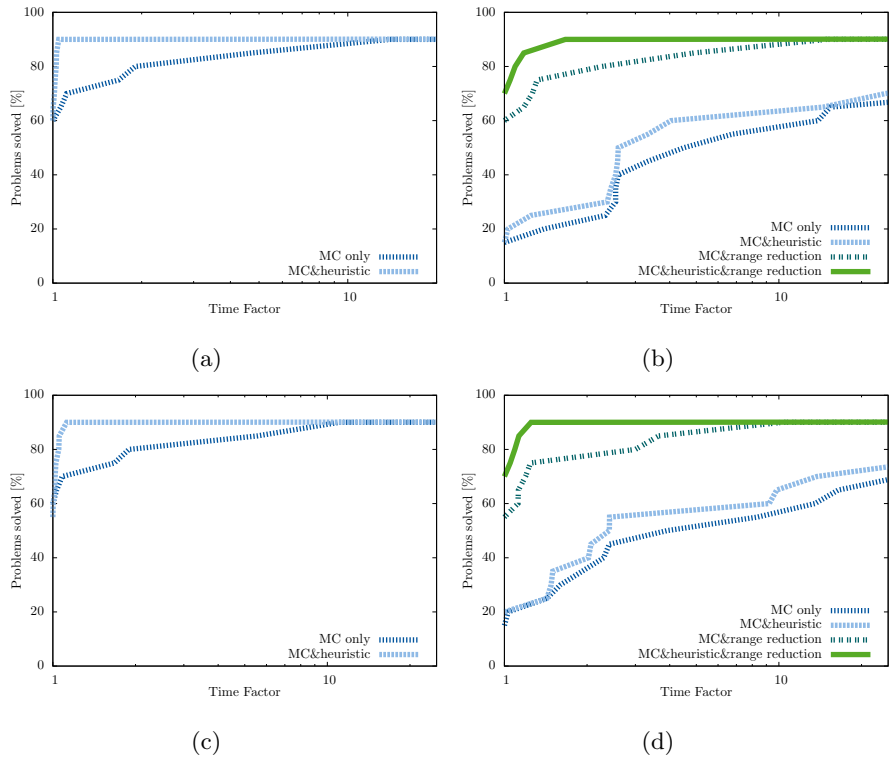

\begin{subfigure}{0.49\textwidth}
\scalebox{0.5}{\input{ perfplot1_1.tex}}
\caption{}\label{Fig:performance_plot_midpoint_a}
\end{subfigure}%
\begin{subfigure}{0.49\textwidth}
\scalebox{0.5}{\input{ perfplot1_2.tex}}
\caption{}\label{Fig:performance_plot_midpoint_b}%
\end{subfigure}
\noindent
\begin{subfigure}{0.49\textwidth}
\scalebox{0.5}{\input{ perfplot1_4.tex}}
\caption{}\label{Fig:performance_plot_incumbent_c}
\end{subfigure}%
\begin{subfigure}{0.49\textwidth}
\scalebox{0.5}{\input{ perfplot1_5.tex}}
\caption{}\label{Fig:performance_plot_incumbent_d}
\end{subfigure}
\caption{\textbf{(a)} Performance plot comparing McCormick relaxation with and without applying the heuristic at the midpoint for $N=1$.\\
\textbf{(b)} Performance plot comparing the additional application of range reduction techniques and the heuristic at the midpoint for $N=1$.\\
\textbf{(c)} Performance plot comparing McCormick relaxation with and without applying the heuristic at the incumbent for $N=1$.\\
\textbf{(d)} Performance plot comparing the additional application of range reduction techniques and the heuristic at the incumbent for $N=1$.}
\label{Fig:performance_plot_midpoint}
\end{figure}

In all numerical experiments we set the absolute and relative optimality tolerances to $\epsilon=10^{-4}$ and absolute and relative feasibility tolerances to $\epsilon=10^{-6}$. First, we compare the impact of the heuristic with only 1 iteration. We allowed for a maximum of 3600 seconds. We compare the number of iterations and the computation time needed for solving the problem when
\begin{itemize}[itemindent=6em]
\item[\textbf{(MC only)}] using McCormick relaxations in a simple B\&B algorithm without the presented heuristic. 
\item[\textbf{(MC heur)}] using McCormick relaxations in a simple B\&B algorithm with the presented heuristic with the middle point of the current node as initial point. 
\item[\textbf{(MC RR)}] using McCormick relaxations in a simple B\&B algorithm. Additionally, also using \emph{Optimization Based Bound Tightening} improved by the filtering bounds technique with factor $0.5$ described in \cite{gleixner2017three} and also employing bound tightening based on the dual multipliers returned by CPLEX\cite{ryoo1995global}. 
\item[\textbf{(MC heur RR)}] using McCormick relaxation with the presented heuristic with the middle point of the current node as initial point. Additionaly, also using Optimization Based Bound Tightening and duality-based bound tightening described in the above point.   
\end{itemize}
Table \ref{table2} in Appendix \ref{appendix} summarizes the results for only 1 allowed iteration within Algorithm \ref{alg:heuristic}. Figures~\ref{Fig:performance_plot_midpoint_a} and \ref{Fig:performance_plot_midpoint_b} show performance plots for the heuristic applied at the midpoint. We observe that the heuristic has the potential to drastically decrease the number of iterations and the computational time needed. Moreover, in the cases where the heuristic did not improve the relaxations, the number of iterations remained the same and the computational time only increased, if at all, by a very marginal amount in all cases. This is explained by the fact that if only 1 iteration of the heuristic is allowed, we can directly integrate the heuristic into the computation of McCormick relaxations and the heuristic only has a constant computational complexity in each propagation step.  

Next, we compare the impact of the heuristic with only 1 iteration but a different initial point. Figures~\ref{Fig:performance_plot_incumbent_c} and \ref{Fig:performance_plot_incumbent_d} show performance plots for the heuristic applied at the incumbent. Again, we allowed for a maximum of 3600 seconds. This time, instead of the middle point of the node, we use the incumbent found by the local solver in the upper bounding procedure as the initial point \emph{and} as the only linearization point in order to construct the linear lower bounding problem. If any coordinate of the incumbent is not within the current node, we simply use the appropriate coordinate of the middle point of this node instead. Table \ref{table3} in Appendix \ref{appendix} summarizes the results for only 1 allowed iteration within Algorithm \ref{alg:heuristic} with the current incumbent as initial point. We see again that the computational time only increased, if at all, by a very marginal amount in all cases. We observe that the choice of the initial point may have a great impact on the advantage provided by the heuristic. Determination of a suitable initial point for the heuristic remains a potential future work and is out of the scope of this article.
\begin{figure}[t]
\centering
\scalebox{0.65}{\input{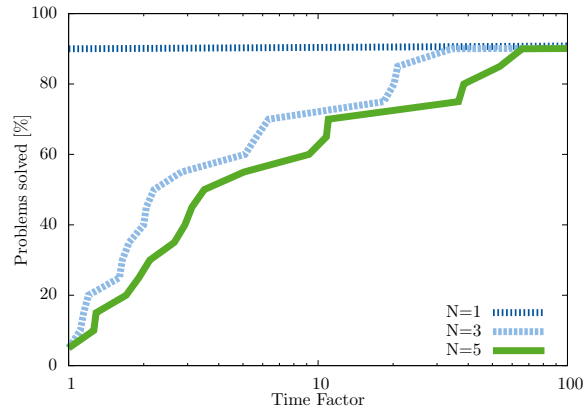}}
\caption{Performance plot comparing the heuristic for $N=1,N=3$ and $N=5$ applied at midpoint without additional range reduction.}
\label{Fig:performance_plot_N>1}
\end{figure}

Last, we compare the impact of the number of iterations within the heuristic. Table \ref{tableIterations} summarizes the numerical results and Figure~\ref{Fig:performance_plot_N>1} shows the performance plot comparing the heuristic for $N=1,N=3$ and $N=5$ applied at midpoint without any additional range reduction.  We allow for a maximum of 3600 seconds. We see that the number of iterations in the B\&B does not increase if we allow more iterations within the heuristic but the computational time needed explodes. This is the behavior that we already shortly discuss in Section \ref{subsec:disc}. Even if the number of factors within the optimization problems is not very high, the additional time needed for further iterations within each factor adds up quickly and is clearly visible if many iterations are needed in order to solve a given optimization problem, see in particular problems \texttt{ex6\_2\_14}, \texttt{growthls}, \texttt{himmelbf} and \texttt{meanvar} in Table \ref{tableIterations} in Appendix \ref{appendix}.

\section{Conclusion}\label{sect:conclusion}

We present a new heuristic for tightening of the univariate McCormick relaxations \cite{mccormick1976,mccormick1983nonlinear} and its extension to multivariate outer functions \cite{tsoukalas2014multivariate} of a factorable function $g$ based on the idea of using tighter interval bounds for the range of each factor of $g$ and obtaining these through subgradients, presented in Section 2.3 in \cite{Najman2017anchoring} and Example 4.4 in \cite{mitsos2009mccormick}. The algorithm  possibly improves the range bounds of the factors of $g$. It uses subgradient propagation for McCormick relaxations \cite{mitsos2009mccormick} in order to construct simple valid affine under- and overestimators of each factor. Then, the affine relaxations are solved with simple function evaluations resulting in  improved range bounds for each factor. This results in tighter McCormick relaxations of the original function $g$. 

Subsequently, we provide numerical results confirming the potential of the presented heuristic. We observe that allowing for only one iteration within the heuristic results in the best computational times. Although more iterations give potentially better bounds, this leads to recalculation of a possibly high number of factors of a factorable function $g$. Moreover, we see that selecting a good initial point may significantly improve the outcome of the algorithm, which remains a potential future work regarding the presented algorithm. Algorithm \ref{alg:heuristic} is especially effective if the underlying interval bounds of the factors are very loose. This is often the case when the well-known dependency problem applies \cite{Moore:1979:MAI:1098639,nla.cat-vn1029192}. A combination of the presented method with the reverse propagation of McCormick relaxations presented in \cite{Wechsung2015} appears to be promising, since it works with range bounds of $g$ and could result in an even greater improvement of McCormick relaxations overall, representing a further potential future development for the McCormick technique. One could also think of a combination of the auxiliary variable method \cite{tawarmalani2002convexification} and the McCormick technique to isolate problematic factors by introduction of auxiliary variables and then directly applying the presented heuristic to these. 

\vspace{0.5cm}
\noindent\textbf{Acknowledgements} We would like to thank Dr.\,Beno\^it Chachuat for providing MC\texttt{++} v2.0 and Dominik Bongartz for providing numerical case studies. This project has received funding from the Deutsche Forschungsgemeinschaft (DFG, German Research Foundation) \textit{Improved McCormick Relaxations for the efficient Global Optimization in the Space of Degrees of Freedom} MA 1188/34-1.

\newpage
\appendix
\section{Appendix}\label{appendix}
\begin{table}[H]
\caption{Table summarizing the problems for the numerical studies for different numbers of iterations within the heuristic.  \#var represents the number of variables, \#ineq stands for the number of inequalities and \#eq for the number of equalities. The domain denotes the domain, we used for the variables, if it was not already \emph{given}. The first 13 problems can be found in the COCONUT benchmark library \cite{Shcherbina2003}.}
\label{table1}
\begin{tabular}{ |l|r|r|r|r|  }
\hline
Name & \#var & \#ineq & \#eq & \multicolumn{1}{c|}{domain} \\
\hline
alkyl     & 14 & 0 & 7 & \emph{given}\\ 
bard      & 3 & 0 & 0 & $[0.001,10]^3$\\
eg1       & 3 & 0 & 0 & $[-10,10]\times[-1,1]\times[1,2]$\\
ex3\_1\_1 & 8 & 6 & 0 & \emph{given}\\
ex6\_2\_10& 6 & 0 & 3 & \emph{given}\\ 
ex6\_2\_14& 4 & 0 & 2 & \emph{given}\\ 
ex8\_1\_3 & 2 & 0 & 0 & $[-2,2]^2$\\
growthls  & 3 & 0 & 0 & $[1,2]\times[0,1]\times[0,1]$\\
himmelbf  & 4 & 0 & 0 & $[0,10]\times[0,10^3]\times[0,10^5]\times[-100,0]$\\
meanvar   & 8 & 0 & 2 & \emph{given}\\
mh4wd     & 5 & 0 & 3 & $[-1000,1000]^5$\\
process   & 10& 0 & 7 & \emph{given}\\
6-hump& 2 & 0 & 0 & $[-3,3]\times[-2,2]$\\
rosenmmx      & 5 & 4 & 0 &$[-10,10]^4\times[-100,100]$\\
ursem\_waves  & 2 & 0 & 0 &$[-0.9,1.2]\times[-1.2,1.2]$\\ 
CS I\cite{bongartz2017deterministic} & 2 & 0 & 9 & \emph{given}\\ 
CS II\cite{bongartz2017deterministic} & 5 & 1 & 12 & \emph{given}\\
CS III\cite{bongartz2017deterministic} & 8 & 22 & 1 & \emph{given}\\
heat \cite{mitsos2009mccormick} & 1 & 0 & 0 & \emph{given}\\ 
kinetic \cite{mitsos2009mccormick} & 3 & 0 & 0 & \emph{given}\\
\hline
\end{tabular}
\end{table}
\begin{sidewaystable}
\centering
\caption{Table summarizing the numerical results. We allowed only 1 iteration within the heuristic and used the middle point of each node as its initial point.  If the time limit of 3600s has been reached, we provide the reached convergence ratio in percent.} \label{table2}
\begin{tabular}{ |l|r|r|r|r|r|r|r|r| }
\hline
$N=1$ & \multicolumn{2}{ c| }{MC only} & \multicolumn{2}{c|}{ MC heur} & \multicolumn{2}{c|}{ MC RR} & \multicolumn{2}{c|}{ MC heur RR} \\
\hline
Name & \#iter & CPU[s] & \#iter & CPU[s] & \#iter & CPU[s] & \#iter & CPU[s]   \\
\hline
alkyl     & 60209 & 24.25& 12645 & 5.13 & 527 & 1.03 & 73   & 0.21 \\\hline 
bard      & 33799 & 5.52 & 33799 & 5.52 & 23 & 0.09 & 23   & 0.09 \\\hline
eg1       & 1039  & 0.1  & 1039  & 0.1  & 71 & 0.03 &  71   & 0.03 \\\hline
ex3\_1\_1 & 19239 & 4.68 & 19239 & 4.68 & 1253 & 1.84 &  1253 & 1.84 \\\hline
ex6\_2\_10& 529289 & 131.15 & 529289 & 137.34  & 101737 & 56.25 &  101737  & 61.54 \\\hline 
ex6\_2\_14& 32537 & 5.85 & 32169 & 5.5  & 1103 & 0.48 & 789  &0.38 \\ \hline
ex8\_1\_3 & 9881  & 1.25 & 993   & 0.09 & 6573 & 1.34 & 459  &  0.15 \\\hline
growthls  & 4385257& 712.03 & 4385257   & 720.6 & 7281 & 3.15 &  7281 & 3.15 \\\hline
himmelbf  & 588859& 68.2 & 588859& 70 & 1425 & 0.63 & 1425   &   0.63 \\\hline
meanvar   & 17381 & 5.21 & 17381 & 5.3  & 1993 & 2.15 & 1993 &  2.15 \\\hline
mh4wd     & 550483& 90.34& 550483& 90.58& 105 & 0.18 & 105  & 0.18 \\\hline
process   & 5751  & 2.57 & 2557  & 1.35 & 379 & 0.75 & 127  & 0.57 \\\hline
rosenmmx      & 1401  & 0.39 & 1401  & 0.39 & 79 & 0.15 & 79   & 0.15 \\\hline
6-hump& 715 & 0.06 & 611 & 0.06 & 281 & 0.06 & 193 & 0.06 \\\hline
ursem\_waves& 97  & 0.1  & 89    & 0.1  & 25 & 0.1 & 23   & 0.1  \\ \hline
CS I~\cite{bongartz2017deterministic} & 629  & 0.39  & 471 & 0.35 & 71 & 0.28 & 57 & 0.28  \\ \hline
CS II~\cite{bongartz2017deterministic} & 355255  & 179.02  & 207979 & 107.12 & 68757 & 60.09 & 28765 & 26.48  \\\hline
CS III~\cite{bongartz2017deterministic} & 2811446  & 81.97\%  & 2613333 & 94.27\% & 1490944 & 85.86\% & 1241856 & 94.84\%  \\\hline
heat~\cite{mitsos2009mccormick} & 129 & 0.34 & 121 & 0.35 & 107 & 0.4 & 103 & 0.4  \\ \hline
kinetic~\cite{mitsos2009mccormick} & 625561 & 96.55\% & 538940 & 97.08\% & 392892 & 96.43\% & 306474 & 97.06\% \\
\hline
\end{tabular}
\end{sidewaystable}
\begin{sidewaystable}
\caption{Table summarizing the numerical results. We allowed only 1 iteration within the heuristic and used the incumbent as its initial point. If the incumbent was not within the node bounds, we used the middle point of the node instead. If the time limit of 3600s has been reached, we provide the reached convergence ratio in percent.} \label{table3}
\centering
\begin{tabular}{ |l|r|r|r|r|r|r|r|r| }
\hline
$N=1$ & \multicolumn{2}{ c| }{MC only} & \multicolumn{2}{c|}{ MC heur} & \multicolumn{2}{c|}{ MC RR} & \multicolumn{2}{c|}{ MC heur RR} \\
\hline
Name & \#iter & CPU[s] & \#iter & CPU[s] & \#iter & CPU[s] & \#iter & CPU[s]   \\
\hline
alkyl     & 49237 & 18.34& 4209  & 1.66 & 257 & 0.54 & 61   & 0.18 \\ \hline%
bard      & 31905 & 4.7  & 31905 & 4.82 & 117 & 0.1 & 117  & 0.12 \\\hline
eg1       & 1035  & 0.1  & 1035  & 0.1  & 143 & 0.07 & 143  & 0.07 \\\hline
ex3\_1\_1 & 16161 & 3.93 & 16161 & 3.93 & 965 & 1.63 & 965 & 1.63  \\\hline
ex6\_2\_10& 653705 & 175.52 & 653705 & 184.46  & 139017 & 75.34 & 139017 & 80.1 \\\hline 
ex6\_2\_14& 32535 & 5.41 & 32143 & 5.5  & 1107 & 0.45 & 813  & 0.4 \\\hline 
ex8\_1\_3 & 10021 & 1 & 1159 & 0.18 & 6605 & 1.23 & 515  & 0.12 \\\hline
growthls  & 4427869& 957.14 & 4427869 & 975.44 & 433 & 0.2 & 433 & 0.25 \\\hline%
himmelbf  & 621203 & 80.7 & 621203 & 90.61 & 4087 & 2.85 & 4087 & 2.85 \\\hline
meanvar   & 8269 & 2.38 & 8269 & 2.51  & 1133 & 1.24 & 1133 & 1.24 \\\hline
mh4wd     & 129795 & 19.5 & 129795 & 19.56 & 689 & 0.46 & 689  & 0.46 \\\hline%
process   & 5069  & 1.76 & 2255  & 0.93 & 437 & 0.56 & 213  & 0.45 \\\hline%
rosenmmx      & 1  & 0.1 & 1  & 0.1 & 1 & 0.1 & 1   & 0.1 \\\hline%
6-hump& 687 & 0.1 & 587 & 0.1 & 289 & 0.1 & 207 & 0.1 \\\hline%
ursem\_waves& 99  & 0.1  & 93    & 0.1  & 27 & 0.1 & 19   & 0.1  \\ \hline%
CS I~\cite{bongartz2017deterministic} & 627  & 0.4  & 429 & 0.37 & 71 & 0.28 & 57 & 0.25  \\ \hline%
CS II~\cite{bongartz2017deterministic} & 474567  & 264.93  & 291759 & 158.12 & 65319 & 58.73 & 15531 & 16.11  \\\hline%
CS III~\cite{bongartz2017deterministic} & 2813962 & 80.17\% & 2842091 & 93.48\% & 1468812 & 88.17\% & 1372106 & 95.6\%  \\\hline%
heat~\cite{mitsos2009mccormick} & 131 & 0.32 & 123 & 0.31 & 109 & 0.37 & 105 & 0.35 \\\hline%
kinetic~\cite{mitsos2009mccormick} & 653067 & 96.57\% & 614659 & 97.09\% & 391701 & 96.43\% & 372522 & 97.08\% \\%
\hline
\end{tabular}
\end{sidewaystable}
\begin{sidewaystable}
\caption{Table summarizing the numerical results for different numbers of iterations within the heuristic. $N$ denotes the number of iterations within the heuristic. \#iter gives the number of iteration needed in the B\&B algorithm. $|\mathbb{F}|$ denotes the number of factors in the given optimization problem.}
\label{tableIterations}
\centering
\begin{tabular}{ |l|r|r|r|r|r|r|r|  }
\hline
 & \multicolumn{2}{ c| }{$N=1$} & \multicolumn{2}{c|}{$N=3$} & \multicolumn{2}{c|}{$N=5$} & \\
\hline
Name & \#iter & CPU[s] & \#iter & CPU[s] & \#iter & CPU[s] & $|\mathbb{F}|$ \\
\hline
alkyl     & 12645 & 5.13 & 12581 & 26.14 & 12577 & 47.2 & 80\\
bard      & 33799 & 5.52 & 33799 &  9.65 & 33799 & 17.23 & 171\\
eg1       & 1039  & 0.1  & 1039  &  0.125& 1039  & 0.17  & 15\\
ex3\_1\_1  &  19239& 4.68 & 19239 &  5.38 & 19239 & 4.212 & 55\\
ex6\_2\_10  & 529289 & 137.34  &  529289 & 302.06 &  529289  & 402.03 & 215\\
ex6\_2\_14& 32169 & 5.5  & 32169 & 114.72 & 32169 & 294.31 & 101\\
ex8\_1\_3 & 993   & 0.09 & 863 & 0.57 & 837 & 0.99 & 51\\ 
growthls  & 4385257 & 720.6 &  4385257 & 1475.94 & 4385257 & 1912.59 & 122\\
himmelbf  & 588859& 70   & 588859 & 1412.64 & 588859 & 2687.12 & 105\\
meanvar   & 17381 & 5.3  & 17381 & 182.38 & 17381 & 348.95 & 164\\
mh4wd     & 550483 & 90.34 & 550483 & 100.87 & 550483 & 114.31 & 37 \\
process   & 2557  & 1.35 & 2545 & 3.82 & 2541 & 6.81 & 73\\
rosenmmx  & 1401  & 0.39 & 1401 & 2.21 &  1401 & 4.21 & 74\\
6-hump& 611 & 0.06 & 611 & 0.12 & 611 & 0.21 & 21\\
ursem\_waves& 89   & 0.1  & 89   & 0.1  & 89   & 0.1 & 25\\
CS I~\cite{bongartz2017deterministic} & 471 & 0.35 & 471 & 6.42  & 471 & 12.8 & 252\\
CS II~\cite{bongartz2017deterministic}  & 207979 & 107.12 & 201697 & 175.89 & 201485 &  226.71 & 333\\
CS III~\cite{bongartz2017deterministic}  & 2613333 & 94.27\% &   1588210   &   93.35\%   &   1333745   &  92.97\%   & 516\\
heat~\cite{mitsos2009mccormick}  & 121 & 0.35 &  117  & 0.56 &  117  &  0.67 & 1201\\
kinetic~\cite{mitsos2009mccormick}  & 538940 & 97.08\% &  209786  &   97.02\%   &   143212   &  96.58\%  & 4516\\
\hline
\end{tabular}
\end{sidewaystable}
\newpage
\bibliographystyle{abbrv}
{\footnotesize
\bibliography{Bibs}

\begin{thebibliography}{10}

\bibitem{androulakis_95_1}
I.~P. Androulakis, C.~D. Maranas, and C.~A. Floudas.
\newblock {$\alpha$BB}: {A} global optimization method for general constrained
  nonconvex problems.
\newblock {\em Journal of Global Optimization}, 7(4):337--363, 1995.

\bibitem{fadbad2016}
C.~Bendtsen and O.~Staunin.
\newblock {\em {F}{A}{D}{B}{A}{D}++, a flexible {C}++ package for automatic
  differentiation. {V}ersion 2.1}, 2012.
\newblock \url{http://www.fadbad.com}. Accessed 18 October 2016.

\bibitem{bertsekas2015convex}
D.~P. Bertsekas.
\newblock {\em Convex optimization algorithms}.
\newblock Athena Scientific Belmont, 2015.

\bibitem{bertsekas2003convex}
D.~P. Bertsekas, A.~Nedi, A.~E. Ozdaglar, et~al.
\newblock {\em Convex analysis and optimization}.
\newblock Athena Scientific, 2003.

\bibitem{bompadre2012convergence}
A.~Bompadre and A.~Mitsos.
\newblock Convergence rate of {M}c{C}ormick relaxations.
\newblock {\em Journal of Global Optimization}, 52(1):1--28, 2012.

\bibitem{convergenceTaylor2013}
A.~Bompadre, A.~Mitsos, and B.~Chachuat.
\newblock Convergence analysis of {T}aylor models and {M}c{C}ormick-{T}aylor
  models.
\newblock {\em Journal of Global Optimization}, 57(1):75--114, 2013.

\bibitem{bongartz2017deterministic}
D.~Bongartz and A.~Mitsos.
\newblock Deterministic global optimization of process flowsheets in a reduced
  space using mccormick relaxations.
\newblock {\em Journal of Global Optimization}, 2017.

\bibitem{brearley1975analysis}
A.~Brearley, G.~Mitra, and H.~P. Williams.
\newblock Analysis of mathematical programming problems prior to applying the
  simplex algorithm.
\newblock {\em Mathematical programming}, 8(1):54--83, 1975.

\bibitem{Brearley1975}
A.~L. Brearley, G.~Mitra, and H.~P. Williams.
\newblock Analysis of mathematical programming problems prior to applying the
  simplex algorithm.
\newblock {\em Mathematical Programming}, 8(1):54--83, 12 1975.

\bibitem{Chachuatmc++}
B.~Chachuat, B.~Houska, R.~Paulen, N.~Peri'c, J.~Rajyaguru, and M.~E.
  Villanueva.
\newblock Set-theoretic approaches in analysis, estimation and control of
  nonlinear systems.
\newblock {\em IFAC-PapersOnLine}, 48(8):981 -- 995, 2015.
\newblock {\url{https://omega-icl.github.io/mcpp/}(Accessed February 2017)}.

\bibitem{comba1993ne}
J.~L.~D. Comba and J.~Stolfi.
\newblock Affine arithmetic and its applications to computer graphics.
\newblock In {\em Proceedings of VI SIBGRAPI (Brazilian Symposium on Computer
  Graphics and Image Processing)}, pages 9--18. Citeseer, 1993.

\bibitem{cornelius1984computing}
H.~Cornelius and R.~Lohner.
\newblock Computing the range of values of real functions with accuracy higher
  than second order.
\newblock {\em Computing}, 33(3-4):331--347, 1984.

\bibitem{de2004affine}
L.~H. De~Figueiredo and J.~Stolfi.
\newblock Affine arithmetic: concepts and applications.
\newblock {\em Numerical Algorithms}, 37(1-4):147--158, 2004.

\bibitem{du1994cluster}
K.~Du and R.~B. Kearfott.
\newblock The cluster problem in multivariate global optimization.
\newblock {\em Journal of Global Optimization}, 5(3):253--265, 1994.

\bibitem{gleixner2017three}
A.~M. Gleixner, T.~Berthold, B.~M{\"u}ller, and S.~Weltge.
\newblock Three enhancements for optimization-based bound tightening.
\newblock {\em Journal of Global Optimization}, 67(4):731--757, 2017.

\bibitem{hamed1993calculation}
A.~S. E.-D. Hamed and G.~P. McCormick.
\newblock Calculation of bounds on variables satisfying nonlinear inequality
  constraints.
\newblock {\em Journal of Global Optimization}, 3(1):25--47, 1993.

\bibitem{hansen1991analytical}
P.~Hansen, B.~Jaumard, and S.-H. Lu.
\newblock An analytical approach to global optimization.
\newblock {\em Mathematical Programming}, 52(1):227--254, 1991.

\bibitem{Cplex2009}
{I}nternational {B}usiness~{M}achines {C}orporation:.
\newblock {\em {I}{B}{M} {I}{L}{O}{G} {C}{P}{L}{E}{X} v12.5. {A}rmonk}, 2009.

\bibitem{Kannan2017}
R.~Kannan and P.~I. Barton.
\newblock The cluster problem in constrained global optimization.
\newblock {\em Journal of Global Optimization}, 5 2017.

\bibitem{Kannan2017bab}
R.~Kannan and P.~I. Barton.
\newblock Convergence-order analysis of branch-and-bound algorithms for
  constrained problems.
\newblock {\em Journal of Global Optimization}, 6 2017.

\bibitem{Kearfott1993}
B.~Kearfott and K.~Du.
\newblock {\em The Cluster Problem in Global Optimization: the Univariate
  Case}, pages 117--127.
\newblock Springer Vienna, Vienna, 1993.

\bibitem{locatelli2013global}
M.~Locatelli and F.~Schoen.
\newblock {\em Global optimization: theory, algorithms, and applications}.
\newblock SIAM, 2013.

\bibitem{maranas1992global}
C.~D. Maranas and C.~A. Floudas.
\newblock A global optimization approach for {L}ennard-{J}ones microclusters.
\newblock {\em The Journal of {C}hemical {P}hysics}, 97(10):7667--7678, 1992.

\bibitem{mccormick1976}
G.~P. McCormick.
\newblock Computability of global solutions to factorable nonconvex programs:
  Part {I}{-}{C}onvex underestimating problems.
\newblock {\em Mathematical Programming}, 10:147--175, 1976.

\bibitem{mccormick1983nonlinear}
G.~P. McCormick.
\newblock {\em Nonlinear programming: Theory, Algorithms, and Applications}.
\newblock Wiley, New York, 1983.

\bibitem{misener2014}
R.~Misener and C.~Floudas.
\newblock Antigone: Algorithms for continuous / integer global optimization of
  nonlinear equations.
\newblock {\em Journal of Global Optimization}, 59(2-3):503--526, 2014.

\bibitem{mitsos2009mccormick}
A.~Mitsos, B.~Chachuat, and P.~I. Barton.
\newblock {M}c{C}ormick-based relaxations of algorithms.
\newblock {\em SIAM Journal on Optimization}, 20(2):573--601, 2009.

\bibitem{Moore:1979:MAI:1098639}
R.~E. Moore and F.~Bierbaum.
\newblock {\em Methods and applications of interval analysis ({S}{I}{A}{M}
  Studies in Applied and Numerical Mathematics)}.
\newblock Society for Industrial \& Applied Math, 1979.

\bibitem{Morrison2016}
D.~Morrison, S.~Jacobson, J.~Sauppe, and E.~Sewell.
\newblock Branch-and-bound algorithms: A survey of recent advances in
  searching, branching, and pruning.
\newblock {\em Discrete Optimization}, 19:79--102, 2 2016.

\bibitem{Najman20161605}
J.~Najman and A.~Mitsos.
\newblock {C}onvergence {O}rder of {M}c{C}ormick {R}elaxations of {LMTD}
  {F}unction in {H}eat {E}xchanger {N}etworks.
\newblock In Z.~Kravanja and M.~Bogataj, editors, {\em 26th European Symposium
  on Computer Aided Process Engineering}, volume~38 of {\em Computer Aided
  Chemical Engineering}, pages 1605 -- 1610. Elsevier, 2016.

\bibitem{Najman2017anchoring}
J.~Najman and A.~Mitsos.
\newblock On {T}ightness and {A}nchoring of {M}c{C}ormick and {O}ther
  {R}elaxations.
\newblock 2017.
\newblock in Revision.

\bibitem{ninin2010reliable}
J.~Ninin, F.~Messine, and P.~Hansen.
\newblock A reliable affine relaxation method for global optimization.
\newblock {\em 4OR}, 13(3):247--277, 2015.

\bibitem{puranik2017bounds}
Y.~Puranik and N.~V. Sahinidis.
\newblock Bounds tightening based on optimality conditions for nonconvex
  box-constrained optimization.
\newblock {\em Journal of Global Optimization}, 67(1-2):59--77, 2017.

\bibitem{puranik2017domain}
Y.~Puranik and N.~V. Sahinidis.
\newblock Domain reduction techniques for global nlp and minlp optimization.
\newblock {\em Constraints}, pages 1--39, 2017.

\bibitem{nla.cat-vn1029192}
H.~Ratschek and J.~Rokne.
\newblock {\em Computer methods for the range of functions}.
\newblock E. Horwood ; Halsted Press Chichester : New York, 1984.

\bibitem{ryoo1995global}
H.~S. Ryoo and N.~V. Sahinidis.
\newblock Global optimization of nonconvex nlps and minlps with applications in
  process design.
\newblock {\em Computers \& Chemical Engineering}, 19(5):551--566, 1995.

\bibitem{ryoo1996branch}
H.~S. Ryoo and N.~V. Sahinidis.
\newblock A branch-and-reduce approach to global optimization.
\newblock {\em Journal of Global Optimization}, 8(2):107--138, 1996.

\bibitem{sahlodin2011Taylor}
A.~Sahlodin and B.~Chachuat.
\newblock Convex/concave relaxations of parametric odes using taylor models.
\newblock {\em Computers {\&} Chemical Engineering}, 35(5):844 -- 857, 2011.
\newblock Selected Papers from ESCAPE-20 (European Symposium of Computer Aided
  Process Engineering - 20), 6-9 June 2010, Ischia, Italy.

\bibitem{schichl2005interval}
H.~Schichl and A.~Neumaier.
\newblock Interval analysis on directed acyclic graphs for global optimization.
\newblock {\em Journal of Global Optimization}, 33(4):541--562, 2005.

\bibitem{scott2011generalized}
J.~K. Scott, M.~D. Stuber, and P.~I. Barton.
\newblock Generalized {M}c{C}ormick relaxations.
\newblock {\em Journal of Global Optimization}, 51(4):569--606, 2011.

\bibitem{Shcherbina2003}
O.~Shcherbina, A.~Neumaier, D.~Sam-Haroud, X.-H. Vu, and T.-V. Nguyen.
\newblock {\em Benchmarking Global Optimization and Constraint Satisfaction
  Codes}, pages 211--222.
\newblock Springer Berlin Heidelberg, Berlin, Heidelberg, 2003.

\bibitem{Shectman1998}
J.~P. Shectman and N.~V. Sahinidis.
\newblock A finite algorithm for global minimization of separable concave
  programs.
\newblock {\em Journal of Global Optimization}, 12(1):1--36, 1 1998.

\bibitem{smith1997global}
E.~M. Smith and C.~C. Pantelides.
\newblock Global optimisation of nonconvex minlps.
\newblock {\em Computers \& Chemical Engineering}, 21:791--796, 1997.

\bibitem{tawarmalani2002convexification}
M.~Tawarmalani and N.~V. Sahinidis.
\newblock {\em Convexification and global optimization in continuous and
  mixed-integer nonlinear programming: theory, algorithms, software, and
  applications}, volume~65.
\newblock Springer Science \& Business Media, 2002.

\bibitem{tawarmalani2005}
M.~Tawarmalani and N.~V. Sahinidis.
\newblock A polyhedral branch-and-cut approach to global optimization.
\newblock {\em Mathematical Programming}, 103(2):225--249, 2005.

\bibitem{tsoukalas2014multivariate}
A.~Tsoukalas and A.~Mitsos.
\newblock Multivariate {M}c{C}ormick {R}elaxations.
\newblock {\em Journal of Global Optimization}, 59:633--662, 2014.

\bibitem{Waechter2006}
A.~W{\"a}chter and L.~T. Biegler.
\newblock On the implementation of an interior-point filter line-search
  algorithm for large-scale nonlinear programming.
\newblock {\em Mathematical Programming}, 106(1):25--57, 3 2006.

\bibitem{wechsung2014}
A.~Wechsung, S.~D. Schaber, and P.~I. Barton.
\newblock The cluster problem revisited.
\newblock {\em Journal of Global Optimization}, 58(3):429--438, 2014.

\bibitem{Wechsung2015}
A.~Wechsung, J.~K. Scott, H.~A.~J. Watson, and P.~I. Barton.
\newblock Reverse propagation of {M}c{C}ormick relaxations.
\newblock {\em Journal of Global Optimization}, 63(1):1--36, 9 2015.

\end{thebibliography}
}
\end{document}